\documentclass[final,1p,times]{elsarticle}
\usepackage{times,amsmath,amsbsy,amssymb}
\usepackage{bm}
\usepackage{color}
\usepackage{mathtools}
\usepackage{diagbox}
\usepackage{subcaption}
\usepackage{graphicx}
\usepackage{epsfig}
\usepackage{epstopdf}
\usepackage{algorithm}
\usepackage{algorithmic}
\usepackage{booktabs}
\usepackage{natbib}
\setcitestyle{comma,sort&compress,numbers,square}
\usepackage{hyperref}
\hypersetup{hypertex=true,
            colorlinks=true,
            linkcolor=blue,
            anchorcolor=blue,
            citecolor=blue}
\usepackage{cleveref}
\usepackage[normalem]{ulem}

\usepackage{geometry}
\usepackage{mathrsfs}

\DeclareMathOperator*{\argmin}{arg\,min}

\everymath{\displaystyle}

\begin{document}
\begin{frontmatter}
\title{A Multigrid Preconditioner for Spatially Adaptive High-order Meshless Method on Fluid-solid Interaction Problems}
\author[WISC]{Zisheng Ye}
\author[TU]{Xiaozhe Hu}
\author[WISC]{Wenxiao Pan\corref{cor}}
\ead{wpan9@wisc.edu}
\cortext[cor]{Corresponding author}
\address[WISC]{Department of Mechanical Engineering, University of Wisconsin-Madison, Madison, WI 53706, USA}
\address[TU]{Department of Mathematics, Tufts University, Medford, MA 02155, USA}

\begin{abstract}
We present a monolithic geometric multigrid preconditioner for solving fluid-solid interaction problems in Stokes limit. The problems are discretized by a spatially adaptive high-order meshless method, the generalized moving least squares (GMLS) with adaptive $h$-refinement. For solving fluid-solid interaction problems, we need to deal with a tightly coupled system consisting of the flow field and solid bodies, resulting in a linear system of equations with a block structure. In Stokes limit, solid kinematics can be dominated by the singularities governing the lubrication effects. Resolving those singularities with adaptive $h$-refinement can lead to an ill-conditioned linear system of equations. The key ingredients of the multigrid preconditioner include the interpolation and restriction operators and the smoothers. For constructing the interpolation and restriction operators, we utilize the geometric information of hierarchical sets of GMLS nodes generated in adaptive $h$-refinement. We build decoupled smoothers through physics-based splitting and then combine them via a multiplicative overlapping Schwarz approach. Through numerical examples with the inclusion of different numbers and shapes of solid bodies, we demonstrate the performance and assess the scalability of the designed preconditioner. As the total degrees of freedom and the number of solid bodies increase, the proposed monolithic geometric multigrid preconditioner can ensure convergence and good scalability when using the Krylov iterative method for solving the linear systems of equations generated from the spatially adaptive GMLS discretization. More specifically, for a fixed number of solid bodies, as the discretization resolution is incrementally refined, the number of iterations of the linear solver can be maintained at the same level, indicating nearly linear scalability of our preconditioner with respect to the total degrees of freedom. When the number of solid bodies $N_s$ increases, the number of iterations is nearly proportional to $\sqrt{N_s}$, implying the sublinear optimality with respect to the number of solid bodies.

\end{abstract}


\end{frontmatter}

\section{Introduction}

In many applications, we need to deal with fluid-solid interactions, where freely moving solid bodies interact with fluid flow through bidirectional hydrodynamic couplings. In Stokes limit, solutions to fluid-solid interaction problems can exhibit singular pressures in the vicinity of sharp corners and in narrow lubrication gaps between solid bodies. If the concentration of solid bodies is dense, solid kinematics can be dominated by lubrication effects. Therefore, solving fluid-solid interaction problems can be challenging especially when there are many solid bodies and the moving solids have complex and feature-rich geometries and display large displacements and rotations. A numerical method must be able to handle the tight coupling between solid kinematics and Stokes flow and to resolve the singularities governing the lubrication effects without invoking \emph{ad hoc} lubrication models. Furthermore, as the number of solid bodies increases, the numerical solver needs to be scalable such that practical applications involving large-scale systems can be accessible in simulations. In our previous work \cite{hu2019spatially}, we have developed a high-order accurate numerical method for solving fluid-solid interaction problems in Stokes limit. The method draws on the generalized moving least squares (GMLS) discretization and adaptive $h$-refinement. The GMLS discretization is built upon a rigorous approximation theory~\cite{wendland2004scattered,mirzaei2012generalized} and by local weighted least square optimization can obtain optimal recovery of target functional over a desired reproduction space directly from a point cloud. Thus, the GMLS discretization does not rely on mesh connectivity between discrete nodes, and hence can circumvent the cost of generating meshes and readily handle the insertion and deletion of nodes in $h$-adaptive refinement. Its polynomial reconstruction property allows high-order accuracy achievable by choosing appropriate polynomial bases. To ensure stable solutions, the div-grad operator in Stokes equations is approximated by a staggered discretization \cite{hu2019spatially,trask2018compatible,trask2017high}; the flexibility of GMLS enables the $\epsilon$-ball graph of neighbor connectivity to be used as a surrogate for primal/dual meshes in, e.g., staggered finite volume methods. To be compatible with the divergence-free constraint for incompressible fluid, GMLS enables the reconstruction of velocity over a space of divergence-free polynomials, which can be efficiently performed locally in a least-squares sense \cite{hu2019spatially,trask2018compatible}, contrary to mesh-based methods necessitating costly global construction of divergence-free shape functions. To achieve adaptive $h$-refinement, an \textit{a posteriori} recovery-based error estimator is defined and employed as the criterion to direct spatial adaptivity \cite{hu2019spatially}. We have demonstrated that the resultant spatially adaptive GMLS method can achieve the same order of accuracy for both the velocity and pressure fields, recover optimal convergence in the presence of singularities, and resolve the singularities governing the lubrication effects without invoking any subgrid-scale lubrication model \cite{hu2019spatially}. Our previous work \cite{hu2019spatially} focused on the accuracy and stability of the solutions to fluid-solid interaction problems and examined benchmark problems with a few solid bodies. In this work, we scale up the spatially adaptive GMLS discretization for solving larger-scale problems and achieve scalability with increasing numbers of solid bodies, while preserving accuracy and stability.  

The linear systems of equations generated from the GMLS discretization need to be solved properly to achieve scalability and efficiency. Krylov methods are usually used for solving large-scale linear systems of equations, for which robust preconditioners must be developed. Since we are dealing with a coupled system consisting of flow field and solid bodies, and the linear system has a block structure, the preconditioner needs to be specially designed. There are mainly two categories of preconditioners for block-structured linear systems: those based on block-factorization approaches (e.g.,~\cite{SilvesterWathen1994,LoghinWathen2004a,SiefertdeSturler2006,MardalWinther2010, MaHuHuXu2016,AdlerGasparHuRodrigoZikatanov2018}) and those based on monolithic multigrid (MG) methods (e.g.,~\cite{BrandtDinar1979,Brandt1984,BraessSarazin1997,Schoberl1998,OosterleeLorenz2006,Janka2008,WangChen2013,ChenHuWangXu2015}). For the GMLS discretization of the Stokes equations, preconditioners based on block factorization were studied in~\cite{trask2018compatible, hu2019spatially}. The block-factorization-based approaches have shown good scalability on pure fluid problems and can handle cases with inclusion of a small number of solid bodies. However, from our in-house numerical experiments, those approaches scale badly with respect to the number of solid bodies. And their convergence become fragile with increasing numbers of solid bodies and levels of refinement. On the other hand, to the best of our knowledge, monolithic multigrid (MG) methods, especially geometric MG (GMG) methods, have not yet been developed for the GMLS discretization, although they have been developed for finite-element methods~\cite{BrandtDinar1979,Brandt1984,BraessSarazin1997,Schoberl1998,OosterleeLorenz2006,Janka2008,WangChen2013,ChenHuWangXu2015}. In a monolithic GMG method, smoothers are usually used for handling all unknowns simultaneously to deal with the block structure of a linear system. The interpolations/restrictions for each unknown are constructed based on a series of hierarchical geometric grids. Such a geometric approach leads to good performance for block-structured problems such as Stokes equations~\cite{BraessSarazin1997,Schoberl1998,OosterleeLorenz2006,WangChen2013,ChenHuWangXu2015} and the fluid-solid interaction problems~\cite{HronTurek2006,TurekHron2010,GeeKuttlerWall2011}.

Due to its meshless nature, there are no hierarchical geometric \emph{meshes} in the GMLS discretization, making the development of a monolithic GMG method challenging. Thus, in our previous works~\cite{trask2018compatible, hu2019spatially}, the algebraic multigrid (AMG) method was used in conjunction with the block-factorization-based preconditioners. In particular, since AMG methods are developed mainly for scalar Poisson-like problems, they are used to invert the diagonal blocks that appear in the block preconditioners. However, due to the complexity of the fluid-solid interaction problems and the limitation of block-factorization-based preconditioners, such AMG-based preconditioning strategies cannot ensure convergence and achieve scalability. Therefore, in this work, we propose to use the geometric information of the adaptively refined GMLS nodes to build the hierarchical structure and develop a monolithic GMG-based preconditioner for solving the linear systems arising from the GMLS discretization.  

The two main ingredients of an MG method are: 1) the interpolation/restriction operators that allow transferring the velocity and pressure values between different resolution levels, and 2) the smoothers that can efficiently damp the high-frequency error components on each resolution level. The interplay between these two ingredients determines the performance of a MG method. We propose a monolithic MG method in the geometric setting in this work. Although we do not have hierarchical meshes, the iterative adaptive $h$-refinement process provides us with hierarchical sets of GMLS nodes whose geometric information can be used to construct the interpolation and restriction operators. More specifically, taking the fine-level nodes as the target nodes and the coarse-level nodes as the source, we use the local GMLS approximations from the source to the target to construct the interpolation operator. In addition, the divergence-free property of velocity, which usually plays an essential role in designing efficient MG methods for solving Stokes equations, is preserved by the GMLS approximations in the interpolation operator. On the other hand, the restriction operator is constructed by an averaging approach, i.e., the value of a (parent) GMLS node with the coarse-level resolution is the averaged value of its child nodes (with one level finer resolution) generated from a new iteration of adaptive refinement. 

For the smoothers, we build decoupled smoothers through physics-based splitting; namely, we handle the fluid domain and the solid bodies separately. For the fluid domain, the corresponding submatrix has a block structure due to the two field variables, velocity and pressure. We use a node-wise block Gauss-Seidel smoother to handle the coupling between the two variables. The smoother for the solid bodies handles each solid body separately. For each solid body, due to the coupling between the solid body and nearby fluids, we assemble the submatrix corresponding to the solid body \emph{and} its neighboring GMLS nodes in the fluid domain. Such a submatrix also has a block structure but is relatively small. So we use a Schur complement approach to invert it approximately based on the block splitting between the fluid and solid degrees of freedoms (DOFs). Finally, the smoother for the fluid domain and the smoother for the solid bodies are combined through a two-stage approach, which can be considered as a multiplicative overlapping Schwarz method.


Built on the proposed interpolation and restriction operators and smoothers, we develop a monolithic GMG preconditioner for using the Krylov iterative method, e.g., GMRES, to solve the linear systems of equations generated from the GMLS discretization. 
With the help of the geometric information of the GMLS nodes, the proposed monolithic GMG preconditioner ensures convergence for complex fluid-solid interaction problems and achieves scalability. Specifically, according to our numerical experiments, for a fixed number of solid bodies, the number of iterations maintains in the same level as we incrementally refine the discretization resolution, which indicates that our preconditioner scales nearly linearly with respect to the number of total DOFs. In addition, when the number of solid bodies $N_s$ increases, the number of iterations is nearly proportional to $\sqrt{N_s}$, which implies the sublinear optimality with respect to the number of solid bodies. Furthermore, the parallel implementation of the proposed monolithic GMG preconditioner is developed by leveraging PETSc~\cite{petsc-user-ref}. A series of numerical tests are designed and performed to demonstrate the parallel scalability of our proposed preconditioner for the inclusion of different shapes of solid bodies. In particular, our parallel implementation achieves weak scalability for both pure fluid flow and fluid-solid interaction problems. Moreover, for the problem in an artificial vascular network, the numerical results are consistent with the strong scalability predicted by Amdahl's law up to $240$ cores, which demonstrates the efficiency of our parallel implementation.

The rest of the paper is organized as follows. \S\ref{sec:preliminary} provides the necessary preliminaries on the governing equations of the fluid-solid interaction problems considered and the spatially adaptive GMLS discretization. In \S\ref{sec:linear_solver}, we present the proposed monolithic GMG method for solving the linear system generated from the GMLS discretization. We describe the parallel implementation  in~\S\ref{sec:implementation_detail}. Our numerical tests and results are discussed in ~\S\ref{sec:numerical_experiment}. Finally, we summarize our main contributions and results and suggest directions of future work in~\S\ref{sec:conclusion}.

\section{Preliminary} \label{sec:preliminary}

\subsection{Governing equations}

We consider the coupled dynamics of steady Stokes flow of incompressible fluid and freely moving solid bodies suspended in the fluid. Hence, the computational domain contains the fluid and $N_s$ solid bodies. 
Denote $\Omega_f$ as the sub-domain occupied by the fluid, $\Gamma_0$ as the outer wall boundary of the whole domain, and $\Gamma_n, n = 1, 2, \dots N_s$ as the boundaries of $N_s$ freely moving solids. Assume $\Gamma_n \cap \Gamma_m = \emptyset, m, n = 0, 1, 2, \dots N_s$ always hold during the whole dynamics process, i.e., there is finite separation between any two boundaries. Each solid body is assumed to follow rigid-body dynamics, and its motion is tracked by its center of mass (COM) position $\mathbf{X}_n$ and orientation $\boldsymbol{\Theta}_n$. Its translational and angular velocities are accordingly donated as $\dot{\mathbf{X}}_n$ and $\dot{\boldsymbol{\Theta}}_n$, respectively. The fluid flow is coupled to the motions of all solid bodies through the following steady Stokes problem:
\begin{equation}
    \left\{
    \begin{aligned}
         & \dfrac{\nabla p}{\rho} - \nu \nabla^2 \mathbf{u} = \mathbf{f}                                                     & \quad & \forall \mathbf{x} \in \Omega_f                    \\
         & \nabla \cdot \mathbf{u} = 0                                                                                       & \quad & \forall \mathbf{x} \in \Omega_f                    \\
         & \mathbf{u} = \mathbf{w}                                                                                           & \quad & \forall \mathbf{x} \in \Gamma_0                   \\
         & \mathbf{u} = \dot{\mathbf{X}}_n + \dot{\boldsymbol{\Theta}}_n \times (\mathbf{x} - \mathbf{X}_n)                  & \quad & \forall \mathbf{x} \in \Gamma_n, n = 1, \dots, N_s \\
         & \mathbf{n} \cdot \dfrac{\nabla p}{\rho} - \nu \mathbf{n} \cdot \nabla ^2 \mathbf{u} = \mathbf{n} \cdot \mathbf{f} & \quad & \forall \mathbf{x} \in \Gamma = \Gamma_0 \cup \Gamma_1 \cup \Gamma_2 \dots \cup \Gamma_{N_s} \;,
    \end{aligned}
    \right.
    \label{eq:governing_eq}
\end{equation}
where $\nu$ is the kinematic viscosity of fluid; $\rho$ is the density of the fluid; $\mathbf{f}$ denotes the body force exerted on the fluid; $\mathbf{w}$ is the velocity of the wall boundary; and 
$\mathbf{n}$ is the unit normal vector outward facing at boundary $\Gamma$.

For a divergence-free velocity field ($\nabla \cdot \mathbf{u}=0$), Eq. \eqref{eq:governing_eq} is equivalent to:
\begin{equation}
    \left\{
    \begin{aligned}
         & \dfrac{\nabla p}{\rho} + \nu \nabla \times \nabla \times \mathbf{u} = \mathbf{f}                                                    & \quad & \forall \mathbf{x} \in \Omega_f                    \\
         & -\dfrac{\nabla^2 p}{\rho} = -\nabla \cdot \mathbf{f}                                                                                & \quad & \forall \mathbf{x} \in \Omega_f                    \\
         & \mathbf{u} = \mathbf{w}                                                                                                             & \quad & \forall \mathbf{x} \in \Gamma_0                    \\
         & \mathbf{u} = \dot{\mathbf{X}}_n + \dot{\boldsymbol{\Theta}}_n \times (\mathbf{x} - \mathbf{X}_n)                                    & \quad & \forall \mathbf{x} \in \Gamma_n, n = 1, \dots, N_s \\
         & \mathbf{n} \cdot \dfrac{\nabla p}{\rho} + \nu \mathbf{n} \cdot \nabla \times \nabla \times \mathbf{u} = \mathbf{n} \cdot \mathbf{f} & \quad & \forall \mathbf{x} \in \Gamma= \Gamma_0 \cup \Gamma_1 \cup \Gamma_2 \dots \cup \Gamma_{N_s} \;,
    \end{aligned}
    \right.
    \label{eq:governing_eq_recast}
\end{equation}
where we have utilized the vector identity $\nabla^2 \mathbf{u} = -\nabla \times \nabla \times \mathbf{u}$ for $\nabla \cdot \mathbf{u}=0$. Solving Eq.~\eqref{eq:governing_eq_recast} instead of Eq. \eqref{eq:governing_eq} can avoid the saddle-point structure of the Stokes operator. However, it necessitates a numerical reconstruction of velocity faithful to the divergence-free constraint.


In Stokes limit, the inertia effect of solids can be neglected. Hence, each solid body's translational and angular motions are subject to the force-free and torque-free constraints as:
\begin{equation}
    \left\{
    \begin{aligned}
         & \int_{\Gamma_n} \boldsymbol{\sigma} \cdot \mathrm{d} \boldsymbol{\mathcal{A}} + \mathbf{f}_{e,n} = \mathbf{0} \\
         & \int_{\Gamma_n} ( \mathbf{x} - \mathbf{X}_n ) \times ( \boldsymbol{\sigma} \cdot \mathrm{d} \boldsymbol{\mathcal{A}}) + \boldsymbol{\tau}_{e,n} = \mathbf{0} \;, \\
    \end{aligned}
    \right.
    \label{eq:freely_solid_dynamics}
\end{equation}
where $n = 1, \dots, N_s$; $\boldsymbol{\sigma} = -p\mathbf{I} + \mu [\nabla \mathbf{u} + {(\nabla \mathbf{u})}^\intercal]$ is the stress exerted by the fluid on the boundary of each solid body; $\mathbf{f}_{e,n}$ and $\boldsymbol{\tau}_{e,n}$ represent the external force and torque applied on each solid body, respectively. 

Therefore, by solving Eqs. \eqref{eq:governing_eq_recast} and \eqref{eq:freely_solid_dynamics} concurrently as a monolithic system, we obtain each solid's translational and angular velocities ${\{} \dot{\mathbf{X}}_n, \dot{\boldsymbol{\Theta}}_n {\}}_{n = 1, 2, \dots, N_s}$ as well as the fluid’s velocity ($\mathbf{u}$) and pressure ($p$) fields at each instant time. Given ${\{} \dot{\mathbf{X}}_n, \dot{\boldsymbol{\Theta}}_n {\}}_{n = 1, 2, \dots, N_s}$ and the initial value of $\{ \mathbf{X}_n, \boldsymbol{\Theta}_n {\}}_{n = 1, 2, \dots, N_s}$, the solid bodies' new positions and orientations ${\{} \mathbf{X}_n, \boldsymbol{\Theta}_n {\}}_{n = 1, 2, \dots, N_s}$ are updated by invoking a temporal integrator. In this work, we adopt the $5^{\text{th}}$ order Runge-Kutta integrator with adaptive time stepping~\cite{dormand1980family}, also known as {\sc ode45} \cite{shampine1997matlab}, and implement it in our solver. To achieve adaptive time stepping, this integrator needs to specify an initial time step and relative error tolerance, which are set as $\Delta t_0 = 0.2s$ and $10^{-5}$, respectively, in this work.

\subsection{Spatially adaptive GMLS discretization}\label{subsec:GMLS}

For numerically solving Eqs. \eqref{eq:governing_eq_recast} and \eqref{eq:freely_solid_dynamics}, we employ the spatially adaptive GMLS method, developed in our previous work~\cite{hu2019spatially}. The linear system resulting from the adaptive GMLS discretization is solved by the proposed scalable linear solver with multigrid preconditioner. Before we introduce the proposed linear solver in \S\ref{sec:linear_solver}, we briefly review the key ingredients of the adaptive GMLS method in this section, including the GMLS approximation, the divergence-free GMLS reconstruction for the velocity field, the staggered discretization of the div-grad operator, and the adaptive refinement scheme. The improvements beyond our previous work \cite{hu2019spatially} are also discussed.

\subsubsection{GMLS approximation}
\label{subsec:GMLS_basics}

The fluid domain $\Omega_f$ is discretized by a set of collocation points (referred to as interior GMLS nodes). All the boundaries $\Gamma_0 \cup \Gamma_1 \cup \Gamma_2 \dots \cup \Gamma_{N_s}$ are discretized by another set of discrete points (referred to as boundary GMLS nodes), where the boundary conditions (BCs) are imposed. For a GMLS node at $\mathbf{x}_i$ and a given scalar function $\psi$ evaluated at its neighbor domain: $\psi_j = \psi(\mathbf{x}_j)$, a polynomial $\psi_h(\mathbf{x})$ of order $m$ is sought to approximate $\psi$ and its derivatives ($D^\alpha\psi$) at $\mathbf{x}_i$. That is, $\psi_h(\mathbf{x}) = \mathbf{P}^\intercal(\mathbf{x}) \mathbf{c}^*$ with a polynomial basis $\mathbf{P}(\mathbf{x})$ and coefficient vector $\mathbf{c}^*$ such that the following weighted residual functional is minimized:
\begin{equation}
    J(\mathbf{x}_i) = \sum_{j \in \mathcal{N}_{\epsilon_i}} {\left[ \psi_j - \mathbf{P}^\intercal_i(\mathbf{x}_j) \mathbf{c}_i \right]}^2 W_{ij} \;.
    \label{eq:GMLS_minization}
\end{equation}
For this quadratic programming optimization problem, its solution can be easily given as:
\begin{equation}
    \mathbf{c}^*_i = {\left( \sum_{k\in\mathcal{N}_{\epsilon_i}} \mathbf{P}_i (\mathbf{x}_k) W_{ik} \mathbf{P}_i^\intercal (\mathbf{x}_k) \right)}^{-1} \left( \sum_{j\in\mathcal{N}_{\epsilon_i}} \mathbf{P}_i^\intercal(\mathbf{x}_j) W_{ij} \psi_j \right) \;.
    \label{eq:c_star_expression}
\end{equation}
Note that for $\psi \in \mathrm{span}(\mathbf{P}(\mathbf{x}))$, this approximation ensures $\psi$ to be exactly reconstructed. This property of polynomial reconstruction grants high-order accuracy achievable for the GMLS approximation by taking large $m$, e.g. $m=4, 6$~\cite{wendland2004scattered}. An arbitrary $\alpha^\text{th}$ order derivative of the function can then be approximated by: \begin{equation}
    D^\alpha \psi(\mathbf{x}) \approx D^\alpha_h \psi_h(\mathbf{x}) \coloneqq {\left(D^\alpha \mathbf{P}(\mathbf{X})\right)}^\intercal \mathbf{c}^*.
    \label{eq:GMLS_derivative}
\end{equation}

The weight function involved in Eqs. \eqref{eq:GMLS_minization}-\eqref{eq:c_star_expression} is defined as $W_{ij} = W(r_{ij})$ with $r_{ij} = \| \mathbf{x}_i - \mathbf{x}_j \|$ and $W(r) = 1 - {\left( \tfrac{r}{\epsilon} \right)}^4$ for $r < \epsilon$, or otherwise, $W(r)= 0$. Here, $\epsilon$ is the compact support, and hence, it is only necessary to include GMLS nodes within an $\epsilon$-neighborhood of the $i$th GMLS node, i.e., $j \in \mathcal{N}_{\epsilon_i} = \{ \mathbf{x}_j \text{ s.t. } \| \mathbf{x}_i - \mathbf{x}_j \| < \epsilon_i \}$. Therefore, the approximation errors in both the interpolant and derivatives of $\psi$ depend on the value of $\epsilon$. On one hand, a small $\epsilon$ leads to high accuracy and sparse linear operators; on the other hand, $\epsilon$ must be large enough to ensure unisolvency over the reconstruction space and thereby a well-posed solution to Eq. \eqref{eq:GMLS_minization}.

\subsubsection{Divergence-free GMLS reconstruction for the velocity field}

We apply the GMLS approximation discussed above to both the velocity and pressure fields. However, the polynomial basis used to approximate the velocity field $\mathbf{u}$ is chosen from the space of $m^\text{th}$ order divergence-free vector polynomials $\mathbf{P}^{\mathrm{div}}(\mathbf{x})$ to enforce compatibility with the divergence-free constraint on the velocity. Thus, the following 
polynomial reconstruction is built for $\mathbf{u}$ and for discretizing its gradient and curl-curl operators:
\begin{equation}
    \begin{aligned}
         & \mathbf{u}_h(\mathbf{x}_i) = {(\mathbf{P}^\mathrm{div}_i)}^\intercal (\mathbf{x}_i) \mathbf{c}_i^{\mathrm{div}*} \quad ~\text{ with } \mathbf{c}_i^{\mathrm{div}*} = \mathbf{M}_i^{\mathrm{div}^{-1}} \left( \sum_{j \in \mathcal{N}_{\epsilon_i}} \mathbf{P}^\mathrm{div}_i (\mathbf{x}_j) \mathbf{u}(\mathbf{x}_j) W_{ij} \right)\;, \\
         & \nabla_h \mathbf{u}_h (\mathbf{x}_i) = {(\nabla \mathbf{P}^\mathrm{div}_i)}^\intercal (\mathbf{x}_i) \mathbf{c}_i^{\mathrm{div}*}\;, \quad {(\nabla \times \nabla \times)}_h \mathbf{u}_h(\mathbf{x}_i) = {(\nabla \times \nabla \times \mathbf{P}^\mathrm{div}_i)}^\intercal (\mathbf{x}_i) \mathbf{c}_i^{\mathrm{div}*} \;,
    \end{aligned}
    \label{eq:GMLS_velocity}
\end{equation}
where
\begin{equation*}
    \mathbf{M}_i^{\mathrm{div}} = \sum_{j \in \mathcal{N}_{\epsilon_i}} \mathbf{P}^\mathrm{div}_i (\mathbf{x}_j) W_{ij} {(\mathbf{P}^\mathrm{div}_i)}^\intercal (\mathbf{x}_j) \;.
\end{equation*}

To impose the Dirichlet (no-slip) BCs in Eq.~\eqref{eq:governing_eq_recast} for $\mathbf{u}$, the boundary GMLS nodes take the velocity of the corresponding boundary $\Gamma_n$ that they belong to.


\subsubsection{Staggered discretization of div-grad operator}\label{sec:staggered_scheme}

To ensure compatibility and thereby numerical stability, we employ the staggered GMLS discretization \cite{trask2017high,hu2019spatially} for the div-grad operator of pressure $p$ in Eq. \eqref{eq:governing_eq_recast}. It builds upon each $\epsilon$-neighborhood graph that plays as a local primal-dual complex for each GMLS node. We refer to it as virtual cell, where a set of primal edges are constructed as: $\mathbf{E}_i = \{ \mathbf{x}_i - \mathbf{x}_j | \mathbf{x}_j \in \mathcal{N}_{\epsilon_i} \}$, each  associated with a midpoint $\mathbf{x}_{ij} = \tfrac{\mathbf{x}_i + \mathbf{x}_j}{2}$ and a virtual dual face at the midpoint and normal to the edge. The staggered GMLS discretization of the div-grad operator is then constructed from a topological gradient over primal edges and a GMLS divergence recovered from local virtual dual faces. Thus, we have:
\begin{equation}
    \begin{aligned}
         & p_h(\mathbf{x}_i) = \mathbf{P}_i^\intercal (\mathbf{x}_i) \mathbf{c}_i^*\;, \quad \quad \nabla_h p_h(\mathbf{x}_i) = \dfrac12 {(\nabla \mathbf{P}_i)}^\intercal (\mathbf{x}_i) \mathbf{c}_i^* \;, \quad \quad \nabla^2_h p_h (\mathbf{x}_i) = \dfrac14 {(\nabla^2 \mathbf{P}_i)}^\intercal (\mathbf{x}_i) \mathbf{c}_i^* \\
        & \text{with} ~\mathbf{c}_i^* = {\left( \sum_{j \in \mathcal{N}_{\epsilon_i}} \mathbf{P}_{i} (\mathbf{x}_{ij}) W_{ij} \mathbf{P}_{i}^\intercal (\mathbf{x}_{ij}) \right)}^{-1} \left( \sum_{j \in \mathcal{N}_{\epsilon_i}} \mathbf{P}_i(\mathbf{x}_{ij}) W_{ij}(p_i - p_j) \right) \;, 
    \end{aligned}
    \label{eq:staggered_operator}
\end{equation}
where the polynomial basis $\mathbf{P}$ can be the $\epsilon_i$-scaled Taylor $m^\text{th}$ order monomials; and  
the coefficient $\mathbf{c}_i^*$ is the solution of the quadratic program:
\begin{equation}\label{eq:quadratic_program_p}
    \mathbf{c}_i(\mathbf{x}) = \argmin_{\mathbf{c}_i} \left\{ \sum_{j \in \mathcal{N}_{\epsilon_i}} {\left[ (p_i - p_j) - \mathbf{P}_i^\intercal(\mathbf{x}_{ij}) \mathbf{c}_i \right]}^2 W_{ij} \right\}.
\end{equation}

Note that pressure $p$ is subject to an inhomogeneous Neumann BC in Eq. \eqref{eq:governing_eq_recast}, written in the form of $\partial_{\mathbf{n}}|_{\Gamma} = g$. To impose this BC, we can add to the quadratic program in Eq. \eqref{eq:quadratic_program_p} for the boundary GMLS nodes the following equality constraint: 
\begin{equation*}
    \mathbf{n}_i \cdot \left[ \dfrac12 {(\nabla \mathbf{P}_i)}^\intercal (\mathbf{x}_i) \mathbf{c}^*_i \right] = g(\mathbf{x}_i)\;, \quad \quad \mathbf{x}_i \in \Gamma \;.
\end{equation*}
In addition, a zero-mean constraint is imposed for $p$ to ensure the uniqueness and physical consistency of the solution. Different from our previous work~\cite{hu2019spatially}, we do not employ a Lagrangian multiplier for enforcing this zero-mean constraint in order to achieve scalability. Details are explained in \S \ref{sec:NBC_implementation}.

\subsubsection{Adaptive refinement}\label{sec:adaptive-algorithm}

In order to reduce computational cost and to recover optimal convergence in the presence of singularities governing lubrication effects, the computational domain and boundaries are discretized with a hierarchy of different resolutions based on an adaptive refinement algorithm \cite{hu2019spatially}.

At each time step, we start from a uniform, coarse initial discretization resolution $\Delta x_i = \Delta x^0$, resulting in a set of GMLS nodes, $\Omega^0$. Then, more GMLS nodes are added adaptively with refined 
$\Delta x_i$, leading to new sets $\Omega^I$ of GMLS nodes. In each refinement iteration, a four-stage procedure is followed:
\begin{equation}
    \mathbf{SOLVE} \to \mathbf{ESTIMATE} \to \mathbf{MARK} \to \mathbf{REFINE} \;.
    \label{eq:adaptive_four_stage}
\end{equation}
In $\mathbf{SOLVE}$ stage, Eq.~\eqref{eq:governing_eq_recast} is numerically solved using the GMLS discretization discussed in \S\ref{subsec:GMLS} and the linear solver introduced in \S\ref{sec:linear_solver}; both the velocity and pressure fields are updated. In $\mathbf{ESTIMATE}$ stage, the recovery-based \emph{a posteriori} error estimator $\eta^r_i$ is evaluated from Eq.~\eqref{eq:local_posteriori_err} for all GMLS nodes. According to the value of the error estimator on each node, we sort all nodes in descending order into a new sequence $\mathbf{x}_{i^\prime}$. In $\mathbf{MARK}$ stage, a fraction of GMLS nodes with largest errors are selected and marked for refining. This fraction of GMLS nodes contribute to $\alpha$ percentage of the total recovered error, where $\alpha$ is a preset parameter and $\alpha=0.8$ in the present work. 
In \textbf{REFINE} stage, the marked GMLS nodes are refined. For solving 2D problems, generally any marked interior GMLS node $\mathbf{x}_i \in \Omega_f$ is split into four nodes with smaller spacing, and a marked boundary node $\mathbf{x}_i \in \Gamma$ is split into two nodes. Denote $\Delta x_i$ as the original spacing (or resolution) of the marked nodes, the nodes newly generated by splitting have the spacing $\Delta x_{i_{new}} = \dfrac12 \Delta x_i$. The updated set of GMLS nodes at the end of \textbf{REFINE} stage is denoted as $\Omega^{I+1}$. These four stages~\eqref{eq:adaptive_four_stage} are repeated iteratively at each time step until the total recovered error (Eq.~\eqref{eq:global_posteriori_err}) is less than a preset tolerance, i.e,
\begin{equation}
    \eta^r \leq \varepsilon_{tol}\;.
    \label{eq:error_tolerance}
\end{equation}


The recovery-based \emph{a posteriori} error estimator is defined based on velocity gradient. Due to the lack of exact solutions in practical applications, true errors cannot be evaluated. Thus, we use the \emph{recovered} error to estimate the true error, which measures the difference between the direct and recovered velocity gradients and can be practically determined in any application problem. More specifically, it is defined as:
\begin{equation}
    \eta_i^r = \dfrac{\sum_{j \in \mathcal{N}_{\epsilon_i}} \| \mathbf{R}{[ \nabla \mathbf{u} ]}_j - \nabla \mathbf{u}_{h, i \to j} \|^2 V_j }{\sum_{j \in \mathcal{N}_{\epsilon_i}} V_j}
    \label{eq:local_posteriori_err} \;.
\end{equation}
Here, the recovered velocity gradient can be evaluated locally on each GMLS node as:
\begin{equation*}
    \mathbf{R}{[\nabla \mathbf{u}]}_i = \dfrac1{N_i} \sum_{j \in \mathcal{N}_{\epsilon_i}} \nabla_h \mathbf{u}_{h, j \to i} \;,
\end{equation*}
where $\nabla_h \mathbf{u}_{h, j \to i} = {(\nabla \mathbf{P}_j^\mathrm{div})}^\intercal (\mathbf{x}_i) \mathbf{c}_j^*$ is the velocity gradient reconstructed at $\mathbf{x}_j$ but evaluated at $\mathbf{x}_i$. In order to properly weigh the contributions from nodes with different discretization resolutions $\Delta \mathbf{x}_i$, a volumetric weight $V_i$ is assigned to each node at $\mathbf{x}_i$ and defined as $V_i = \Delta x_i^d$, where $d$ is the dimension of the problem's physical space. Then, the total recovered error over all GMLS nodes is:
\begin{equation}
    \eta^r = \dfrac{\sum_{\mathbf{x}_i \in \Omega} \eta_i^r V_i}{\sum_{\mathbf{x}_i \in \Omega} \| \nabla \mathbf{u}_{h, i} \|^2 V_i} \;.
    \label{eq:global_posteriori_err}
\end{equation}

For applications of dilute suspensions of solids, the above adaptive refinement procedure can start from an uniform, coarse initial discretization. However, for concentrate suspensions of solids, or when gaps between some solid boundaries are rather narrow, even smaller than $\Delta x^0$, 
the initial coarse discretization can result in none GMLS nodes allocated within the gaps. In that case, adaptive refinement would not take place within those gaps. Thus, we invoke a preprocessing step to guarantee GMLS nodes allocated everywhere throughout the entire computational domain including narrow gaps, before we conduct adaptive refinement following the algorithm discussed above. 

In this preprocessing step, we examine the $\epsilon$-neighborhood of each boundary GMLS node and check if the following requirement is satisfied:
\begin{equation}\label{eq:autoref_req1}
    \text{For}~\mathbf{x}_i \in \Gamma_{n}~\text{and}~\forall j \in \mathcal{N}_{\epsilon_i}, ~~\mathbf{x}_j \notin \Gamma_m ~~ \text{with}~~m\neq n~\text{and}~m, n = 0, 1, \dots, N_s \;.
\end{equation}
This requirement ensures that each boundary GMLS node has enough interior GMLS nodes in its $\epsilon$-neighborhood. If it is not satisfied, the corresponding boundary node and all its neighbor nodes are refined. Note that in this preprocessing step, the neighbor nodes of a boundary node can be within solid bodies.

After preprocessing and any iteration of refinement, we perform a post-processing step to enforce quasi-uniform discretization within any $\epsilon$-neighborhood $\mathcal{N}_{\epsilon_i}$, because large difference in the discretization resolution within an $\epsilon$-neighborhood would result in ill-conditioned GMLS approximation. Thus, for any $\epsilon$-neighborhood that does not satisfy: 
\begin{equation}\label{eq:autoref_req2}
    \dfrac{\max(\Delta x_j)}{\min(\Delta x_k)} \leq 2\;,~~ \text{for}~~\forall \mathbf{x}_i \in \Omega_f\cup \Gamma ~~ \text{and}~~ j, k \in \mathcal{N}_{\epsilon_i} \;,
\end{equation}
we will mark and refine the nodes with the coarsest resolution in that $\epsilon$-neighborhood, until the requirement in Eq. \eqref{eq:autoref_req2} is satisfied.

\subsubsection{Numerical quadrature}

With the GMLS discretization, a composite quadrature rule is employed for approximating the integrals in Eq. \eqref{eq:freely_solid_dynamics} as:
\begin{equation}\label{eq:int_discretization}
    \begin{aligned}
        \int_{\Gamma_n} \boldsymbol{\sigma} \cdot \mathrm{d} \boldsymbol{\mathcal{A}} \approx                                      & \sum_{\mathbf{x}_i \in \Gamma_n} (-p_i \mathbf{I} + \nu [\nabla_h \mathbf{u}_{h, i} + {(\nabla_h \mathbf{u}_{h, i})}^\intercal]) \cdot (\Delta \mathcal{A}_i \mathbf{n}_i)                          \\
        \int_{\Gamma_n} (\mathbf{x} - \mathbf{X}_n) \times (\boldsymbol{\sigma} \cdot \mathrm{d} \boldsymbol{\mathcal{A}}) \approx & \sum_{\mathbf{x}_i \in \Gamma_n} \left( (\mathbf{x}_i - \mathbf{X}_n) \times (-p_i \mathbf{I} + \nu [\nabla_h \mathbf{u}_{h, i} + {(\nabla_h \mathbf{u}_{h, i})}^\intercal]) \right) \cdot (\Delta \mathcal{A}_i\mathbf{n}_i) \;,
    \end{aligned}
\end{equation}
where $\nabla_h \mathbf{u}_{h, i} = {(\nabla \mathbf{P}^\mathrm{div}_i)}^\intercal (\mathbf{x}_i) \mathbf{c}_i^*$; and $\Delta\mathcal{A}_i=\Delta x_i$ for 2D problems. As shown in~\cite{hu2019spatially}, this quadrature rule is sufficient for the entire numerical method to achieve high-order convergence. 

\section{Linear Solver with Multigrid Preconditioner}\label{sec:linear_solver}

Solving the linear systems resulting from discretization dominates the entire computational cost in simulations. Therefore, the main focus of this work is to design a scalable preconditioner for the Krylov method to solve the linear systems generated in the \textbf{SOLVE} stage. A simple preconditioner, such as the Gauss-Seidel preconditioner, cannot ensure convergence and is not effective in practice. In the previous work~\cite{trask2018compatible, hu2019spatially}, a block-factorization-based preconditioner based on decoupling of the velocity and pressure fields was employed, and AMG methods were applied for inverting each block. For solving benchmark smaller-scale fluid-solid interaction problems \cite{trask2018compatible, hu2019spatially}, such a linear solver can work reasonably well. However, as the number of solid bodies increases, e.g., in the examples considered in the present work, the performance of such a block preconditioner deteriorates. One major reason is that the velocity block used in the block-factorization-based preconditioner consists of unknowns related to both the fluid and solid bodies. In particular, it combines the discretized curl-curl operator in Eq.~\eqref{eq:governing_eq_recast} for the fluid with the discretized integrals in  Eq.~\eqref{eq:int_discretization} for the solid bodies. In addition, the inclusion of many solid bodies results in a very irregular shape for the overall computational domain. All of these make it difficult for an AMG method to find a proper coarsening. Consequently, it fails to converge for the velocity block and strongly affects the overall performance of the block-factorization-based preconditioner. 

This section introduces a monolithic GMG method to build a preconditioner for the Krylov linear solver such as GMRES to solve the linear systems resulted from GMLS discretization. The new multigrid method utilizes the hierarchical sets of GMLS nodes generated during the adaptive refinement. Therefore, no coarsening is needed, avoiding the major difficulty that the AMG method previously encountered. Using the geometric information of these hierarchical sets of the GMLS nodes, we can construct interpolation operators based on the GMLS approximation and restriction operators based on averaging. Furthermore, a two-stage smoother is developed based on the splitting between the fluid domain and the solid bodies. Finally, the interpolation/restriction operators and the smoother are used in a V-cycle fashion to define the overall monolithic GMG preconditioner. 

We first discuss about the block structure of the resultant linear system in \S\ref{subsection:structure_linear_system}. Then, we briefly review the required operations and entire process of multigrid preconditioning in \S\ref{sec:multigrid_review}. In \S\ref{subsubsec:restriction_interpolation_operators}, we explain how to construct the interpolation and restriction operators. In \S\ref{subsec:smoother}, we present the smoothers designed based on physics splitting. The entire monolithic GMG method introduced in this work is finally summarized in \S\ref{subsec:preconditioner_summary}.


\subsection{Block structure of the linear system}\label{subsection:structure_linear_system}

After the governing equations~\eqref{eq:governing_eq_recast}-\eqref{eq:freely_solid_dynamics} are discretized by the GMLS discretization discussed in \S\ref{subsec:GMLS}, the resulting linear system has the following block structure:
\begin{equation}
    \mathbf{A} \boldsymbol{\chi} =\mathbf{y}
    \quad \text{with}~~
    \mathbf{A} =
    \begin{bmatrix}
        \mathbf{K} & \mathbf{G} & \mathbf{C} \\
        \mathbf{B} & \mathbf{L}              \\
        \mathbf{D} & \mathbf{T}
    \end{bmatrix}\;, \quad
    \boldsymbol{\chi}=
    \begin{bmatrix}
        \mathbf{u} \\
        p          \\
        \dot{\mathbf{X}}
    \end{bmatrix}\;, \quad
    \mathbf{y} = \begin{bmatrix}
        \mathbf{f}_{tot} \\
        g                \\
        \mathbf{f}_{s}
    \end{bmatrix}
    \;.
    \label{eq:structure_linear_system}
\end{equation}
Here, $\dot{\mathbf{X}}$ includes all solids' both translational and rotational velocities, i.e., $\dot{\mathbf{X}} = {\{ \dot{\mathbf{X}}_n, \dot{\boldsymbol{\Theta}}_n \}}_{n = 1, 2, \dots, N_s}$. $\mathbf{K}$ corresponds to the discretized curl-curl operator ($\nabla \times \nabla \times$) of velocity in Eq.~\eqref{eq:governing_eq_recast}. It is noted that the curl-curl operator for a divergence-free polynomial basis is equivalent to Laplacian. $\mathbf{L}$ corresponds to the discretized Laplacian operator ($\nabla^2$) of pressure in Eq.~\eqref{eq:governing_eq_recast}, which is obtained from the staggered discretization of div-grad operator discussed in~\S\ref{sec:staggered_scheme}. As such, the nonzero diagonal blocks in $\mathbf{A}$, i.e., $\mathbf{K}$ and $\mathbf{L}$, are all discretized Laplacian operators. While $\mathbf{B}$ denotes the contribution from the $\nu \nabla \times \nabla \times$ operator in the inhomogeneous Neumann BC in Eq.~\eqref{eq:governing_eq_recast}, $\mathbf{G}$ represents the discretized $\tfrac1\rho \nabla$ operator. For the discretized force-free and torque-free constraints as in Eq.~\eqref{eq:int_discretization}, $\mathbf{D}$ denotes the contribution from viscous stress with velocity gradient; $\mathbf{T}$ denotes the contribution from pressure. $\mathbf{C}$ contains the velocity constraints (no-slip BCs) from each solid's kinematics on their boundaries in Eq.~\eqref{eq:governing_eq_recast}.  $\mathbf{f}_{tot}$ combines the body force $\mathbf{f}$ exerted in fluid and the velocity $\mathbf{w}$ on the wall boundary. $g$ contains $\nabla \cdot \mathbf{f}$ and $\mathbf{n} \cdot \mathbf{f}$ in Eq.~\eqref{eq:governing_eq_recast}. Finally, $\mathbf{f}_{s}$ represents the external force and torque applied on the solid bodies, such that $\mathbf{f}_{s} = {[-\mathbf{f}_{e, n} \  -\boldsymbol{\tau}_{e,n}]}^\intercal_{n = 1, 2, \dots, N_s}$.

From its block structure, we see that $\mathbf{A}$ is neither symmetric nor positive definite. The strong coupling between the fluid and solid DOFs deteriorates the conditioning of $\mathbf{A}$. Thus, the linear system in Eq. \eqref{eq:structure_linear_system} is challenging to solve in practice.    Furthermore, since enough neighbor nodes are needed in GMLS discretization, the matrix $\mathbf{A}$ is dense. Therefore, designing a robust and efficient linear solver is crucial for achieving scalability in practical simulations.


\subsection{Multigrid methods}
\label{sec:multigrid_review}

Most relaxation-type iterative solvers, like Gauss-Seidel (GS), converge slowly for linear systems with fine resolution discretized from partial differential equations, e.g., \eqref{eq:governing_eq_recast}. The main reason is that the convergence rates of different error components can vary significantly, if we decompose the error using the eigenvectors of the linear system, e.g., $\mathbf{A}$ in Eq.~\eqref{eq:structure_linear_system}. The components with the eigenvectors corresponding to large eigenvalues (in terms of magnitude) are usually referred to as high-frequency components of the error; the components with the eigenvectors corresponding to small eigenvalues are called low-frequency components of the error. In general, for relaxation-type iterative solvers, the high-frequency components converge to zero quickly independent of the discretization resolution. However, the low-frequency components converge slowly, getting worse when the resolutions are refined and hence slowing down the overall convergence. On the other hand, low-frequency components of the error on fine resolutions can be well approximated on coarse discretization and hence become high-frequency error components on coarse resolutions~\cite{Brandt1980}. This fact motivates the idea of moving to a coarser resolution to eliminate the low-frequency error components and leads to a multilevel approach known as the multigrid method \cite{Brandt1980,Brandt1984}.

Multigrid methods exploit a discretization with different resolutions to obtain optimal convergence rate and hence are naturally compatible with adaptive $h$-refinement discussed in \S\ref{sec:adaptive-algorithm}. The operations of going back and forth between coarse and fine resolutions are called the restriction and interpolation, respectively. To minimize the approximation errors across different resolution levels, a smoothing procedure is usually executed before the restriction operation, and another smoothing step is applied after the interpolation operation. The entire process as such is referred to as ``V-cycle"~\cite{Brandt1984} and is summarized in \Cref{algorithm:standard_multigrid}. In the present work, we follow this V-cycle for building our multigrid preconditioner.

\begin{algorithm}
    \caption{V-cycle process of multigrid preconditioning: V-Cycle($\mathbf{A}^I, \mathbf{y}^I$)}\label{algorithm:standard_multigrid}
    \hspace*{\algorithmicindent} \textbf{Input}: Coefficient matrix $\mathbf{A}^I$ and right-hand side (RHS) vector $\mathbf{y}^I$ \\
    \hspace*{\algorithmicindent} \textbf{Output}: Approximate solution $\boldsymbol{\chi}^I$
    \begin{algorithmic}[1]
        \IF{$I == 0$}
        \STATE Solve: $\mathbf{A}^I \boldsymbol{\chi}^I = \mathbf{y}^I$
        \ELSE{}
        \STATE Pre-smooth: $\boldsymbol{\chi}^I = \text{Smooth}(\mathbf{A}^I, \mathbf{y}^I)$
        \STATE Compute residual: $\mathbf{r}^I = \mathbf{y}^I - \mathbf{A}^I \boldsymbol{\chi}^I$
        \STATE Restrict: $\mathbf{y}^{I-1} = \boldsymbol{\mathcal{R}}^I \mathbf{r}^I$
        \STATE Apply V-cycle recursively: $\mathbf{z}^{I-1} = \text{V-Cycle}(\mathbf{A}^{I-1}, \mathbf{y}^{I-1})$
        \STATE Interpolate: $\boldsymbol{\chi}^I = \boldsymbol{\chi}^I + \boldsymbol{\mathcal{I}}^{I} \mathbf{z}^{I-1}$
        \STATE Post-smooth: $\boldsymbol{\chi}^I = \boldsymbol{\chi}^I + \text{Smooth}(\mathbf{A}^I, \mathbf{y}^I - \mathbf{A}^I \boldsymbol{\chi}^I)$
        \ENDIF{}
    \end{algorithmic}
    \hspace*{\algorithmicindent} \textbf{Return}: $\boldsymbol{\chi}^I$
\end{algorithm}

In the monolithic setting, the interpolation operators $\boldsymbol{\mathcal{I}}^I$, the restriction operators $\boldsymbol{\mathcal{R}}^I$, and the smoothers (i.e., subroutine Smooth($\mathbf{A}$, $\mathbf{y}$)) need to be carefully designed to properly handle the coupling of fluid and solid-body DOFs. The adaptive refinement discussed in~\S\ref{sec:adaptive-algorithm} yields a hierarchical series of node sets $\Omega^I$, $I=0,1,2,\dots$. We can utilize the hierarchical structures between node sets to construct the corresponding interpolation $\boldsymbol{\mathcal{I}}^I$ and restriction $\boldsymbol{\mathcal{R}}^I$ operators for each node set $\Omega^I$. 

\subsection{Interpolation and restriction operators}\label{subsubsec:restriction_interpolation_operators}

\begin{figure}[htp]
    \begin{subfigure}{.45\textwidth}
        \centering
        \includegraphics[width=.98\textwidth]{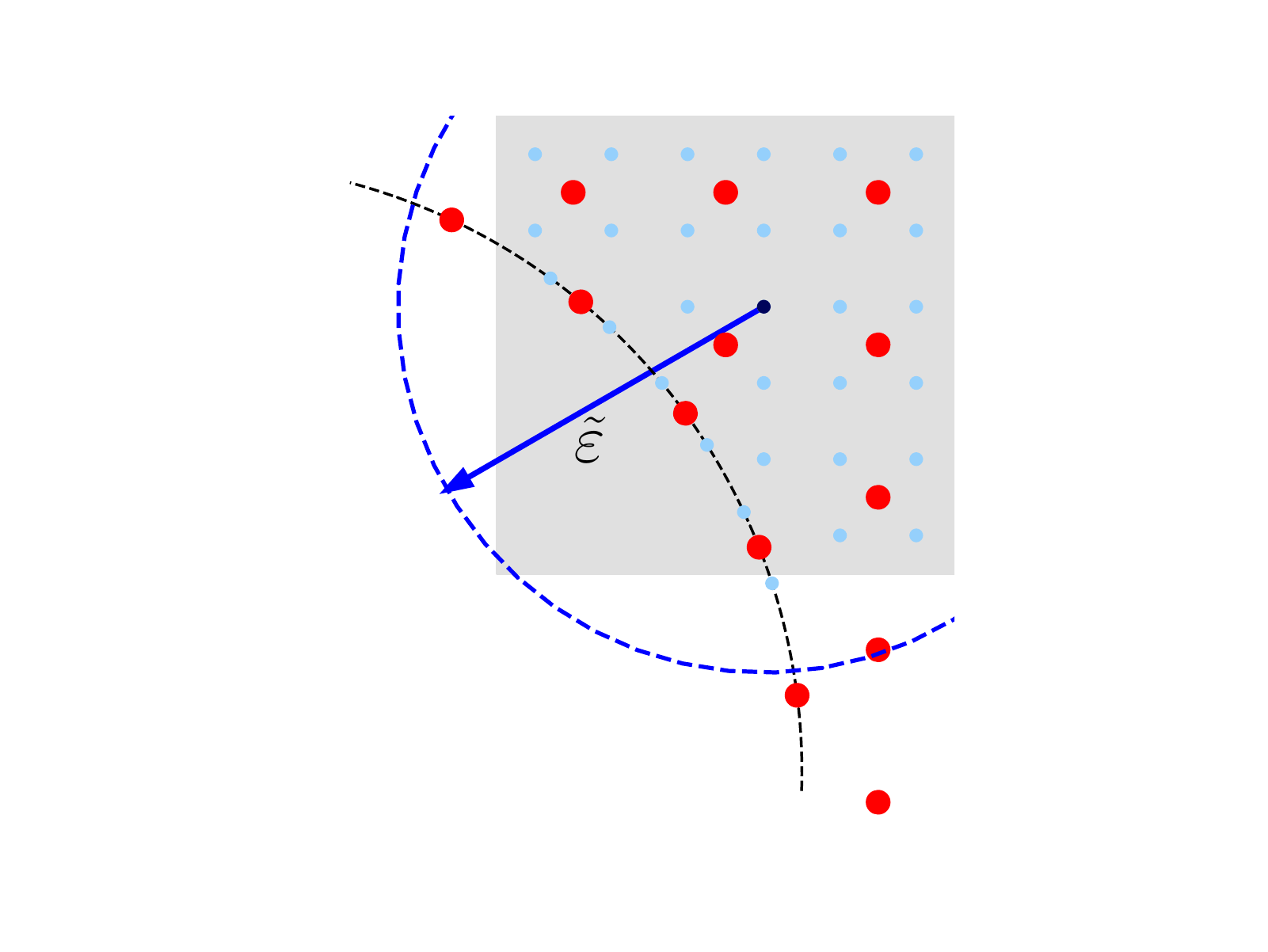}
        \caption{The interpolation operator at a fine-level node $\mathbf{x}_i^F$ (highlighted in dark blue) is constructed as the GMLS approximation from the coarse-level nodes (red) within the fine-level node's $\tilde{\epsilon}$-neighborhood.}
    \end{subfigure}
    \hfill
    \begin{subfigure}{.45\textwidth}
        \centering
        \includegraphics[width=.98\textwidth]{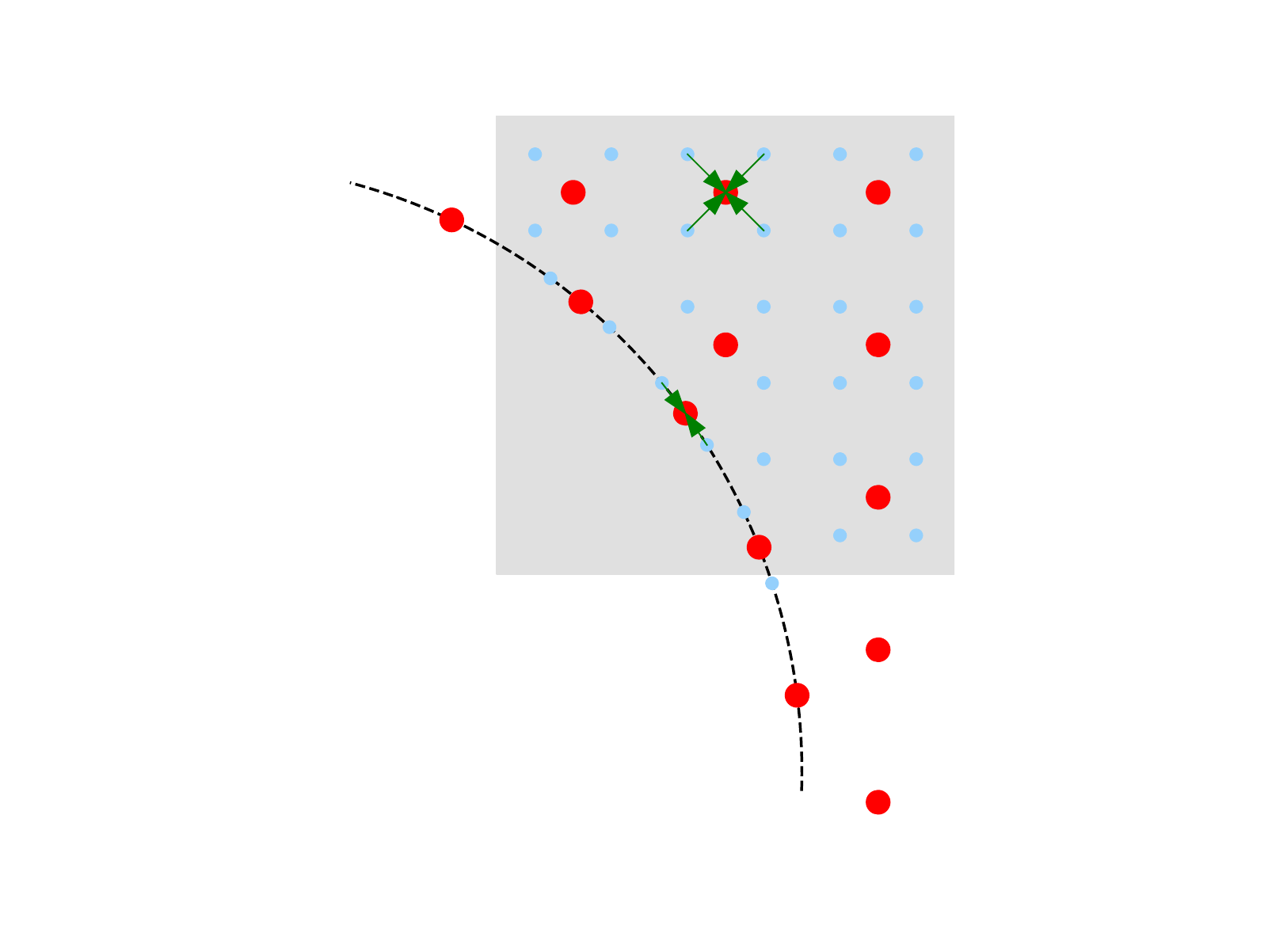}
        \caption{The restriction operator at a coarse-level node $\mathbf{x}_i^C$ (red) is constructed as averaging over its own child nodes (light blue with green arrows) generated in \textbf{REFINE} stage during an iteration of adaptive refinement.}
    \end{subfigure}
    \caption{Schematic of the interpolation and restriction. Coarse-level nodes are displayed in red and fine-level nodes are in blue. The black dashed lines indicate part of a solid body's boundary. }\label{fig:interpolation_restriction_operator}
\end{figure}

In mesh-based discretization methods, e.g. finite element method, the interpolation/restriction operators can be constructed from the nested function spaces across different levels of meshes. In this work, instead of using function spaces, we build the interpolation and restriction operators through local approximations, i.e., approximating any scalar or vector field locally from a set of neighbor nodes. 
More specifically, if we take the fine-level nodes $\Omega^{F}$ as the ``target" nodes and the coarse-level nodes $\Omega^{C}$ as the ``source", by constructing GMLS approximation (\S\ref{subsec:GMLS_basics}) on the target from the source, we define the interpolation operator. If the target are coarse-level nodes, and the source are fine-level nodes, then a local averaging approximation on the target from the source builds the restriction operator. An illustrative description about the interpolation and restriction is provided in ~\Cref{fig:interpolation_restriction_operator}. 


Thus, the interpolation operator for the pressure field is constructed as:
\begin{equation*}
    p (\mathbf{x}_i^F) = \mathbf{P}^\intercal_i(\mathbf{x}_i^F) \mathbf{c}_i^*~~ \text{ with }~ \mathbf{c}_i^* = \mathbf{M}^{-1}_i \left( \sum_{j \in \mathcal{N}_{\tilde{\epsilon}_i}} \mathbf{P}_i (\mathbf{x}_j^C) p(\mathbf{x}_j^C) W_{ij} \right),~~    \mathbf{M}_i = \sum_{j \in \mathcal{N}_{\tilde{\epsilon}_i}} \mathbf{P}_i (\mathbf{x}_j^C) W_{ij} \mathbf{P}^\intercal_i (\mathbf{x}_j^C),~~ \mathbf{x}_i^F \in \Omega^F,~~ \mathbf{x}_j^C \in \Omega^C\;,
\end{equation*}

where the polynomial basis $\mathbf{P}$ is the $\tilde{\epsilon}$-scaled Taylor monomials. 
For velocity, the interpolation operator is built by directly using the divergence-free reconstruction space as follows:
\begin{equation*}
    \mathbf{u}(\mathbf{x}_i^F) = {(\mathbf{P}^{\mathrm{div}}_i)}^\intercal (\mathbf{x}_i^F) \mathbf{c}_i^{\mathrm{div}*} ~~ \text{ with }~ \mathbf{c}_i^{\mathrm{div}*} = \mathbf{M}_i^{\mathrm{div}^{-1}} \left( \sum_{j \in \mathcal{N}_{\tilde{\epsilon}_i}} \mathbf{P}^{\mathrm{div}}_i (\mathbf{x}_j^C) \mathbf{u}(\mathbf{x}_j^C) W_{ij} \right), ~~\mathbf{M}_i^{\mathrm{div}} = \sum_{j \in \mathcal{N}_{\tilde{\epsilon}_i}} \mathbf{P}^{\mathrm{div}}_i (\mathbf{x}_j^C) W_{ij} {(\mathbf{P}^\mathrm{div}_i)}^\intercal (\mathbf{x}_j^C), ~~\mathbf{x}_i^F \in \Omega^F,~~ \mathbf{x}_j^C \in \Omega^C\;.
\end{equation*}
As such, we ensure the divergence-free constraint for velocity is consistently preserved during interpolation across different fine/coarse levels of nodes. As is well-known~\cite{Schoberl1998}, for solving the incompressible Stokes equations, interpolations that maintain the divergence-free property are crucial for developing an efficient GMG preconditioner. In the above GMLS approximations for constructing the interpolation operator, $\tilde{\epsilon}_i$ for each fine-level node needs to be large enough such that sufficient coarse-level nodes are included in each $\tilde{\epsilon}$-neighborhood to ensure unisolvency over the reconstruction space and a well-posed solution
to the weighted least square optimization. 

For the restriction operator, 
since the matrix $\mathbf{A}$ in the linear system ~\eqref{eq:structure_linear_system} is non-symmetric, it is unnecessary to require the restriction operator to be the transpose of the interpolation operator. Hence, we construct the restriction operation based on $h$-refinement. A coarse-level node $\mathbf{x}_i^C \in \Omega^C$ marked in the \textbf{MARK} stage of the adaptive refinement algorithm is called a ``parent" node. The newly generated nodes (within the fluid domain or on solid boundaries) in the \textbf{REFINE} stage are called the ``child nodes" corresponding to their parent node $\mathbf{x}_i^C$. As a result, an interior parent node generally has four child nodes, and a boundary parent node has two child nodes in 2D. The approximation at a parent node is given by averaging the field values from its child nodes, which provides the restriction operator. By such, the matrix corresponding to the restriction operator is 
much sparser than that of the interpolation operator, leading to cheaper cost for matrix multiplication during the restriction operations.

For the blocks in $\mathbf{A}$ related to solid DOFs, the interpolation and restriction operators are simply identity matrices. Assembled together, the interpolation and restriction operators can be summarized as:
\begin{equation}
    \begin{bmatrix}
        \mathbf{u}^{I+1} \\
        p^{I+1}          \\
        \dot{\mathbf{X}}
    \end{bmatrix}
    =
    \boldsymbol{\mathcal{I}}^I
    \begin{bmatrix}
        \mathbf{u}^I \\
        p^I          \\
        \dot{\mathbf{X}}
    \end{bmatrix}
    =
    \begin{bmatrix}
        \boldsymbol{\mathcal{I}}_{\mathbf{u}}                 \\
         & \boldsymbol{\mathcal{I}}_p              \\
         &              & \mathbf{I}
    \end{bmatrix}
    \begin{bmatrix}
        \mathbf{u}^I \\
        p^I          \\
        \dot{\mathbf{X}}
    \end{bmatrix}
    \;,
    \quad \quad
    \begin{bmatrix}
        \mathbf{u}^I \\
        p^I          \\
        \dot{\mathbf{X}}
    \end{bmatrix}
    =
    \boldsymbol{\mathcal{R}}^I
    \begin{bmatrix}
        \mathbf{u}^{I+1} \\
        p^{I+1}          \\
        \dot{\mathbf{X}}
    \end{bmatrix}
    =
    \begin{bmatrix}
        \boldsymbol{\mathcal{R}}_\mathbf{u}               \\
         & \boldsymbol{\mathcal{R}}_p              \\
         &            & \mathbf{I}
    \end{bmatrix}
    \begin{bmatrix}
        \mathbf{u}^{I+1} \\
        p^{I+1}          \\
        \dot{\mathbf{X}}
    \end{bmatrix} \;,
    \label{eq:interpolation_restriction_operator}
\end{equation}
where the super index $I$ corresponds to the node set $\Omega^I$ resulted from the $I$-th iteration of adaptive refinement; while $\boldsymbol{\mathcal{I}}_\mathbf{u}$ denotes the interpolation for velocity, $\boldsymbol{\mathcal{I}}_p$ represents the interpolation for the pressure field; $\mathbf{I}$ is an identity matrix; $\boldsymbol{\mathcal{R}}_p$ is the restriction operator for a scalar valued function (e.g., pressure); $\boldsymbol{\mathcal{R}}_\mathbf{u} = \operatorname{diag}(\boldsymbol{\mathcal{R}}_p, \cdots, \boldsymbol{\mathcal{R}}_p)$ contains $d$ copies of $\boldsymbol{\mathcal{R}}_p$ on the diagonal and is the restriction operator for velocity (a $d$-dimensional vector field).

\subsection{Smoother based on physics splitting}\label{subsec:smoother}

In this subsection, we introduce the design of the smoother, i.e., the subroutine Smooth($\mathbf{A}$, $\mathbf{y}$) used in Algorithm~\ref{algorithm:standard_multigrid}. As we discussed, it is a relaxation-type iterative method that can efficiently smooth the high-frequency components of the error. For solving a fluid-solid interaction problem, we need the smoother to handle the high-frequency components of the errors for the fluid part and solid part, respectively, as well as the strong coupling between them.

Eqs. \eqref{eq:governing_eq_recast}-\eqref{eq:freely_solid_dynamics} inherently state two types of physics: One is Stokesian flow in a confined space, and the other describes the dynamics of several rigid solid bodies undergoing external loads as well as drags exerted by surrounding Stokesian flow. The coupling of different types of physics inspires the design of our physics-based smoother, which contains two stages. The first stage takes care of the fluid DOFs, and the second handles the solid bodies as well as their neighboring fluid nodes. The details are presented as follows.

\subsubsection{Smoother for the fluid DOFs}
The submatrix
\begin{equation}
    \mathbf{F} = \boldsymbol{\mathcal{B}}_F \mathbf{A} \boldsymbol{\mathcal{B}}^\intercal_F
    =
    \begin{bmatrix}
        \mathbf{K} & \mathbf{G} \\
        \mathbf{B} & \mathbf{L}
    \end{bmatrix}
    \label{eq:fluid_part}
\end{equation}
is considered as the fluid part of the system.  Here, $\boldsymbol{\mathcal{B}}_F$ is a Boolean matrix, which picks the fluid field variables, $\mathbf{u}$ and $p$, out of the whole unknown vector. $\mathbf{F}$ contains all the field variables (velocity and pressure) directly related to the interior GMLS nodes. The diagonal part $\mathbf{K}$ and $\mathbf{L}$ denotes the curl-curl operator onto the velocity field and the Laplacian operator onto the pressure field, respectively. As mentioned in~\S\ref{subsection:structure_linear_system}, the curl-curl operator $\mathbf{K}$ is equivalent to the Laplacian operator on the divergence-free polynomial basis. Therefore, the diagonal blocks of the fluid part matrix $\mathbf{F}$ are all Laplacian-type operators. Further, as seen in the governing equation~\eqref{eq:governing_eq_recast}, there exists a strong coupling between the velocity and pressure fields. We suggest using a node-wise Gauss-Seidel (GS) smoother. In 2D, we organize all field values into a $3\times 3$ sub-matrix for each GMLS node. Therefore, we would include its diagonal entities related to the two Laplacian operators and the discretized coupling terms at the node in the smoother. From our numerical experiments, the node-wise smoother outperforms the normal entry-wise one, and inverting small dense matrices does not affect the scalability of the smoother. This node-wise GS smoother is denoted as $\mathbf{S}_{\mathbf{F}}$ in \Cref{algorithm:smoother}.
 
\subsubsection{Smoother for the solid bodies}

\begin{figure}[htp]
    \centering
    \includegraphics[width=8cm]{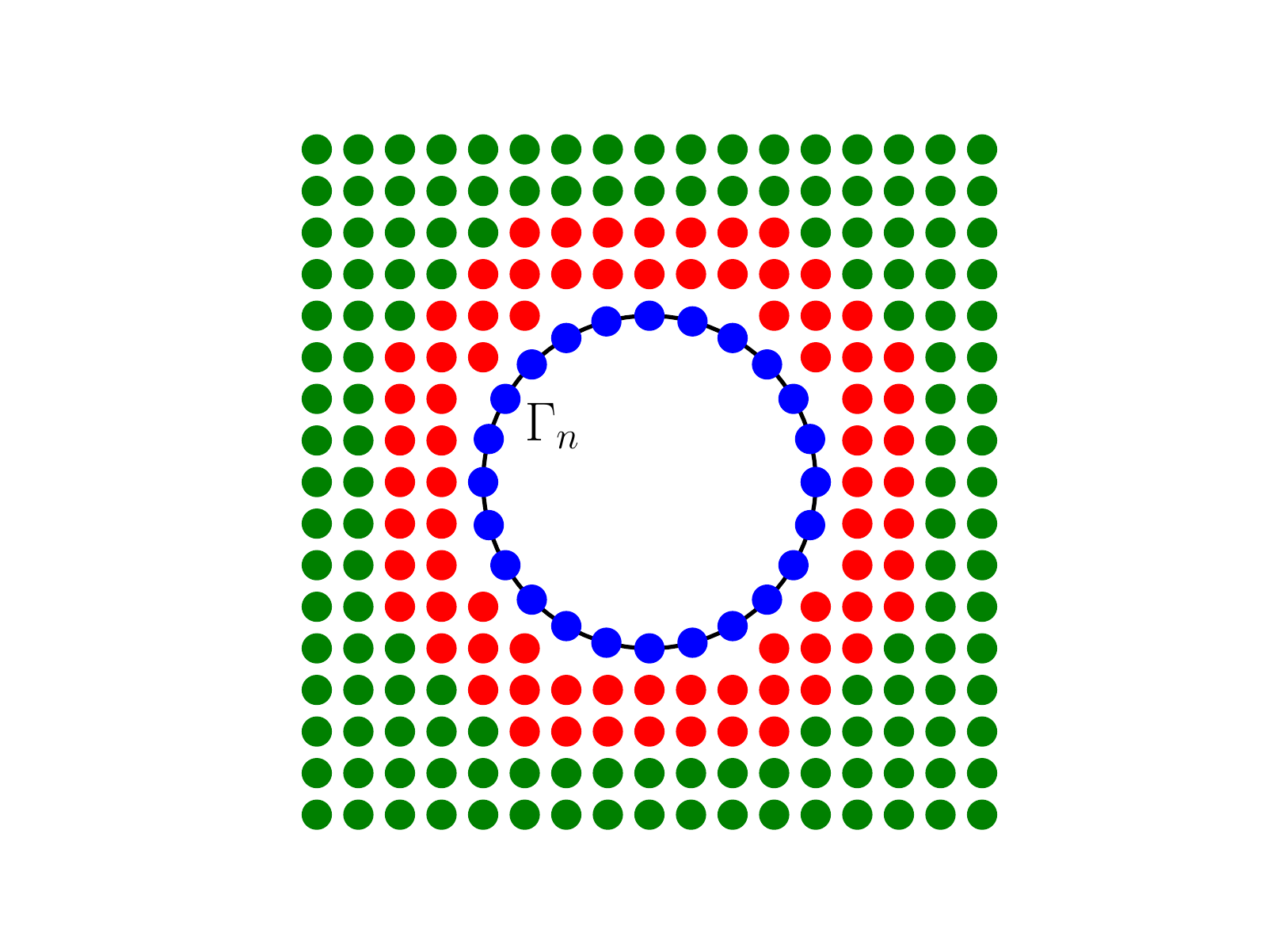}
    \caption{Illustration of construction of $Q_n$.  Blue nodes denote the boundary nodes on $\Gamma_n$.  Red nodes represent the nodes near the boundary $\Gamma_n$ and contribute to the force and torque terms to the $n$-th solid.  The green nodes denote the normal interior GMLS nodes.}
    \label{fig:solid_neighbor}.
\end{figure}

Now consider the sub-matrix directly related to each solid body and its neighbor fluid nodes. The submatrix
\begin{equation}
    \mathbf{N}_n =
    \boldsymbol{\mathcal{B}}_{N_n} \mathbf{A} \boldsymbol{\mathcal{B}}^\intercal_{N_n}
    =
    \begin{bmatrix}
        \mathbf{K}_n & \mathbf{G}_n & \mathbf{C}_n \\
        \mathbf{B}_n & \mathbf{L}_n                \\
        \mathbf{D}_n & \mathbf{T}_n
    \end{bmatrix}
    \label{eq:neighboring_matrix}
\end{equation}
corresponds to the $n$-th solid body and its neighbor fluid nodes. It is a square matrix and contains discretized drag force and torque exerted by the fluid. Each $\mathbf{N}_n$ is constructed by the following procedure. First, for each solid body, we construct an index set $Q_n$ such that
\begin{equation*}
    Q_n = \{i, j~|~\mathbf{x}_i \in \Gamma_n,~ j \in \mathcal{N}_{\epsilon_i} \}\;, \quad n = 1, 2, \dots, N_s \;.
\end{equation*}
\Cref{fig:solid_neighbor} illustrates the construction of $Q_n$, where the green and red nodes are the interior GMLS (fluid) nodes; the blue nodes are the GMLS nodes on $\Gamma_n$ (the boundary of the $n$-th solid body).  
Since the red nodes are within the $\epsilon$-neighborhood of the blue nodes, $Q_n$ includes the indices of the nodes rendered in blue and red in \Cref{fig:solid_neighbor} but excludes the nodes rendered in green. Second, we form the Boolean matrix $\boldsymbol{\mathcal{B}}_{N_n}$ from $Q_n$, which picks the variables related to the $n$-th solid body and its neighbor fluid nodes out of the whole unknown vector. The submatrix $\mathbf{N}_n$ is then built according to Eq. \eqref{eq:neighboring_matrix}. Note that there is no intersection between most $Q_n$ for $n = 1, 2, \dots, N_s$, especially when after several iterations of adaptive refinement, there are plenty of neighbor fluid nodes between any two closely contacting solid bodies. Thus, the parallel scalability can still be ensured. 

Once we construct all sub-matrices $\mathbf{N}_n$, we use a Schwarz-type relaxation method to build the smoother. The resulting additive Schwarz-type smoother can be defined as $\mathbf{S}_{\mathbf{N}}^{\text{exact}}:= \sum_{n=1}^{N_s} \mathcal{B}_{N_n}^{\intercal} \mathbf{N}_n^{-1} \mathcal{B}_{N_n}$. Exactly inverting $\mathbf{N}_n$ could be expensive. Therefore, we use a Schur complement approach to invert it approximately, i.e., 
\begin{equation}
	\widetilde{\mathbf{N}}_n^{-1} =
	\begin{bmatrix}
		\mathbf{I} & -\widetilde{\mathbf{F}}_n^{C^{-1}}
		\begin{bmatrix}
			\mathbf{C}_n \\
			\mathbf{0}
		\end{bmatrix} \\
		\mathbf{0} & \mathbf{I}
	\end{bmatrix}
	\begin{bmatrix}
		\widetilde{\mathbf{F}}_n^{C^{-1}} & \mathbf{0} \\
		\mathbf{0} & \widetilde{\boldsymbol{\Psi}}_n^{-1}
	\end{bmatrix}
	\begin{bmatrix}
		\mathbf{I} & \mathbf{0} \\
		-
		\begin{bmatrix}
			\mathbf{D}_n & \mathbf{T}_n
		\end{bmatrix}
		\widetilde{\mathbf{F}}_n^{C^{-1}}
		& \mathbf{I}
	\end{bmatrix}\;,
	\label{eq:schur_preconditioner}
\end{equation}
where $\widetilde{\mathbf{F}}_n^{C^{-1}}$ approximates the inverse of the submatrix
\begin{equation*}
	\mathbf{F}_n^C = 
	\begin{bmatrix}
		\mathbf{K}_n & \mathbf{G}_n \\
		\mathbf{B}_n & \mathbf{L}_n
	\end{bmatrix}\;.
\end{equation*}
In our numerical experiments, we use one iteration of the node-wise GS smoother to define $\widetilde{\mathbf{F}}_n^{C^{-1}}$. In addition, the approximated Schur complement $\widetilde{\boldsymbol{\Psi}}_n$ is given as:
\begin{equation*}
	\widetilde{\boldsymbol{\Psi}}_n =
	\begin{bmatrix}
		\mathbf{D}_n & \mathbf{T}_n
	\end{bmatrix}
	\text{block-diag}^{-1}(\mathbf{F}_n^C)
	\begin{bmatrix}
		\mathbf{C}_n \\
		\mathbf{0}
	\end{bmatrix}\;.
\end{equation*}
Here, $\text{block-diag}^{-1}(\cdot)$ denotes the diagonal block inversion of a matrix, which takes a $3\times3$ sub-matrix for each GMLS node in $\mathbf{F}_n^C$ and inverts those sub-matrices. Since the size of matrix $\widetilde{\boldsymbol{\Psi}}_n$ is quite small, only $3 \times 3$ in 2D, a direct solver is applied whenever the inversion of $\widetilde{\boldsymbol{\Psi}}_n$ is needed.  Finally, for the solid bodies, the overall additive Schwarz-type smoother using the approximate Schur complement approach is given by:
\begin{equation}\label{eq:add-schwarz-solid}
	\mathbf{S}_{\mathbf{N}} := \sum_{n=1}^{N_s} \mathcal{B}_{N_n}^{\intercal} \widetilde{\mathbf{N}}_n^{-1} \mathcal{B}_{N_n} \;.
\end{equation}

Given the two smoothers $\mathbf{S}_{\mathbf{F}}$ and $\mathbf{S}_{\mathbf{N}}$, we connect them through an overlapping multiplicative Schwarz approach and hence establish the proposed two-stage smoothing scheme. \Cref{algorithm:smoother} summarizes the established overall smoother. 
\begin{algorithm}[H]
    \caption{Smoother: $\text{Smooth}(\mathbf{A}, \mathbf{y})$}\label{algorithm:smoother}
    \hspace*{\algorithmicindent} \textbf{Input}: Coefficient matrix $\mathbf{A}$ and RHS vector $\mathbf{y}$ \\
    \hspace*{\algorithmicindent} \textbf{Output}: Relaxed solution $\boldsymbol{\chi}$
    \begin{algorithmic}[1]
    	\STATE Initialize $\boldsymbol{\chi} = \mathbf{0}$
        \STATE Perform $k$ steps of node-wise GS smoother on the fluid DOFs:
        \FOR{$i=1,\cdots,k$} 
        \STATE $\boldsymbol{\chi} \gets \boldsymbol{\chi} + \mathcal{B}_F^{\intercal}  \mathbf{S}_{\mathbf{F}} \mathcal{B}_{F} (\mathbf{y} - \mathbf{A} \boldsymbol{\chi}) $
        \ENDFOR
        \STATE Apply the additive Schwarz-type smoother as in Eq.~\eqref{eq:add-schwarz-solid} for the solid bodies:
        \vskip -10pt 
        $$\boldsymbol{\chi} \gets \boldsymbol{\chi} + \mathbf{S}_{\mathbf{N}} (\mathbf{y} - \mathbf{A} \boldsymbol{\chi}), \quad 
        \text{where} \ \mathbf{S}_{\mathbf{N}} = \sum_{n=1}^{N_s} \mathcal{B}_{N_n}^{\intercal} \widetilde{\mathbf{N}}_n^{-1} \mathcal{B}_{N_n}
        $$
        \STATE \vskip -10pt  \textbf{Return}: $\boldsymbol{\chi}$
    \end{algorithmic}
\end{algorithm}

\subsection{Linear solver summary}\label{subsec:preconditioner_summary}

Given the constructed interpolation/restriction operators (Eq. \eqref{eq:interpolation_restriction_operator}) and smoothers (\Cref{algorithm:smoother}), we follow the V-cycle process in \Cref{algorithm:standard_multigrid} to build the entire monolithic GMG preconditioner. Since the coefficient matrix $\mathbf{A}$ is non-symmetric, we employ the GMRES method as the Krylov iterative solver for solving Eq. \eqref{eq:structure_linear_system}. In each GMRES iteration, we call once the multigrid preconditioning, i.e.,  \Cref{algorithm:standard_multigrid}. With the proposed monolithic GMG preconditioner, we aim to ensure the convergence of the linear solver and to optimize the scaling of the number of GMRES iterations required with respect to the numbers of solid bodies and total DOFs.

\section{Parallel Implementation}\label{sec:implementation_detail}
To solve large-scale fluid-solid interaction problems, the parallel implementation of the proposed monolithic GMG preconditioner is developed. We aim to achieve the parallel scalability of our numerical solver.



\subsection{Domain decomposition and neighbor search}
All the GMLS nodes, $\forall \mathbf{x}_i \in \Omega^I$, yield after each iteration $I$ of adaptive refinement, are evenly distributed into $N_c$ sets, where $N_c$ is the number of CPU cores invoked in the simulation. Each set of nodes, denoted as $\Omega_k^I$ ($k=1, 2, \dots, N_c$), is then allocated into one core, following the Recursive Coordinate Bisection (RCB) method~\cite{zoltan2-website}. In each core, $\mathtt{Compadre}$ \cite{paul_kuberry_2019_3338664} is used to generate the coefficients resulted from the GMLS discretization for the node-set $\Omega_k^I$. This way of domain decomposition guarantees that the nodes in the same set $\Omega^I$ are evenly distributed among CPU cores and spatially clustered on each core, and hence balances the workload and minimizes the communication required between cores for solving the linear system in Eq. \eqref{eq:structure_linear_system}.

After domain decomposition, a ghost node set $\mathcal{G}_{k_\xi}^I$ is determined according to the order ($m$) of GMLS discretization and the spatial locations of the GMLS nodes in the node set $\Omega_k^I$, as:
\begin{equation*}
    \mathcal{G}_{k_\xi}^I = \{\mathbf{x}_j \ | \ \|\mathbf{x}_j - \mathbf{x}_i\| < \xi \Delta x_i, ~\forall \, \mathbf{x}_i \in \Omega_k^I, ~\mathbf{x}_j \in \Omega^I  \}\;,~~k=1, 2, \dots, N_c\; ,~~ \xi = \xi(m) \;,
\end{equation*}
where $\xi$ is a function of $m$ and we use $\xi = 4m$ in our numerical tests; $\Delta x_i$ denotes the discretization resolution of node $\mathbf{x}_i$. Collecting this ghost node set only calls one communication between all cores. Once the collection is done, neighbor search for each node in $\Omega_k^I$ is performed by calling  $\mathtt{nanoflann}$~\cite{blanco2014nanoflann}, which is a function building KD-trees~\cite{bentley1975multidimensional} and finds all nodes in the $\epsilon$-neighborhood of any node $\mathbf{x}_i \in \Omega_k^I$ from the ghost node set $\mathcal{G}_{k_\xi}^I$, i.e., 
\begin{equation*}
    \mathcal{N}_{\epsilon_i} = \{ j \ | \ \| \mathbf{x}_i - \mathbf{x}_j \| < \epsilon_i, ~\forall \, \mathbf{x}_i \in \Omega_k^I, ~ \mathbf{x}_j \in \mathcal{G}_{k_{\xi}}^I \} \;.
\end{equation*}
$\mathcal{N}_{\epsilon_i}$ hence provides the neighbor list needed in the GMLS discretization that leads to the coefficient matrix $\mathbf{A}$ of the linear system in Eq. \eqref{eq:structure_linear_system}.   


Note that when building the interpolation operator in \S\ref{subsubsec:restriction_interpolation_operators}, GMLS approximation is needed, for which we also need to build the neighbor list. Different from the above, the neighbor nodes of $\mathbf{x}_i$ in this GMLS approximation are not in the same node set as $\mathbf{x}_i$. Thus, for that purpose the ghost node set is determined as: 
\begin{equation*}
    \widetilde{\mathcal{G}}_{k_\xi}^I = \{\mathbf{x}_j \ | \ \|\mathbf{x}_j - \mathbf{x}_i\| < \xi \Delta x_i, ~\forall \mathbf{x}_i \in \Omega_k^{I+1}, ~\mathbf{x}_j \in \Omega^I \} \;,  ~~k=1, 2, \dots, N_c\; ,
\end{equation*}
where $\mathbf{x}_i$ and $\mathbf{x}_j$ belong to the node sets yield from different iterations of adaptive refinement, respectively. And then the neighbor list is built. This time, $\mathtt{nanoflann}$ finds all nodes in the $\epsilon$-neighborhood of any node $\mathbf{x}_i \in \Omega_k^{I+1}$ in the ghost node set $\widetilde{\mathcal{G}}_{k_\xi}^I$. Thus, we have: 
\begin{equation*}
    \widetilde{\mathcal{N}}_{\epsilon_i} = \{ j \ | \ \| \mathbf{x}_i - \mathbf{x}_j \| < \epsilon_i, ~\forall \mathbf{x}_i \in \Omega_k^{I+1}, ~ \mathbf{x}_j \in \widetilde{\mathcal{G}}_{k_\xi}^I \} \;,
\end{equation*}
which provides the neighbor list for the GMLS approximation needed to construct the interpolation operator. 



\subsection{Data storage}\label{subsec:data_storage}
To reduce the data storage, all DOFs related to the solid bodies can be stored in a single CPU core (e.g., the last core). Thus, following the convention of PETSc, the sub-matrices $\mathbf{D}$ and $\mathbf{T}$ reside in the last rows of the matrix $\mathbf{A}$. However, the number of nonzero entities in $\mathbf{D}$ and $\mathbf{T}$ would increase with the inclusion of more solid bodies. As a result, the parallel scalability deteriorates as the number of solid bodies increases. Therefore, we spread the whole storage of data and workload related to $\mathbf{D}$ and $\mathbf{T}$ among all cores as follows. When assembling $\mathbf{D}$ and $\mathbf{T}$, each core would go through the GMLS nodes stored in it to generate its local $\mathbf{D}_i$ and $\mathbf{T}_i$. Accordingly, any matrix-vector multiplication operation for $\mathbf{D}$ and $\mathbf{T}$ would be split among all cores. That is, the matrix-vector multiplication is done locally within each core, followed by a parallel summation over all cores. For example, the matrix-vector multiplication for $\mathbf{D}$ is given by:
\begin{equation}
    \mathbf{D} \mathbf{u} = \sum_{i = 0}^{N_c - 1} \mathbf{D}_i \mathbf{u} \;,
    \label{eq:rb_mat}
\end{equation}
where $N_c$ is the total number of CPU cores invoked in a simulation; $\mathbf{D}_i \mathbf{u}$ is computed in each core parallelly and then summed up over all cores. By such, the parallel scalability is recovered as the number of solid bodies increases.

\subsection{Linear algebra operations}
The parallel implementation of our linear solver is achieved by interfacing with PETSc package~\cite{petsc-user-ref}. Once all coefficients resulting from the GMLS discretization are generated from each core, the entire coefficient matrix $\mathbf{A}$ is assembled in parallel through PETSc. And all linear algebra operations associated with the proposed monolithic GMG preconditioner, as well as the Krylov iterations, are performed by calling the corresponding PETSc functions, including the matrix-vector multiplication, matrix-matrix multiplication, and matrix inversion by a direct solver, as well as the node-wise GS iteration and the GMRES iterative solver.

\subsection{Neumann BC for the pressure} \label{sec:NBC_implementation}
As stated in Eq.~\eqref{eq:governing_eq_recast}, an inhomogeneous Neumann BC is imposed for pressure. 
To ensure the uniqueness and physical consistency of the solution, we additionally enforce a zero-mean constraint for the pressure field, i.e., requiring $\boldsymbol{\Xi}^\intercal \mathbf{p} = 0$, where $\boldsymbol{\Xi}$ is a constant vector (e.g. $\boldsymbol{\Xi} = \mathbf{1}$); $\mathbf{p}$ denotes the vector of discretized pressure with entities $p(\mathbf{x}_i)$, $\mathbf{x}_i \in \Omega$. If we use a Lagrange multiplier to impose this zero-mean constraint, it would break the block structure of $\mathbf{F}$ in Eq. \eqref{eq:fluid_part}, and in the meanwhile introduce a dense row and column, i.e., $\boldsymbol{\Xi}^\intercal$ and $\boldsymbol{\Xi}$, respectively, into $\mathbf{F}$, which in turn would deteriorate the parallel scalability of our preconditioner. Thus, we instead solve the following problem:
\begin{equation}
    \left(\mathbf{I} - \dfrac{\boldsymbol{\Xi} \boldsymbol{\Xi}^\intercal}{\boldsymbol{\Xi}^\intercal \boldsymbol{\Xi}} \right) \mathbf{L} \mathbf{p} = \left(\mathbf{I} - \dfrac{\boldsymbol{\Xi} \boldsymbol{\Xi}^\intercal}{\boldsymbol{\Xi}^\intercal \boldsymbol{\Xi}} \right) \mathbf{g} \;,
    \label{eq:projected_system}
\end{equation}
which is equivalent to using the Lagrangian multiplier. Here, $\mathbf{g}$ denotes the vector with entities $g(\mathbf{x}_i)$, $\mathbf{x}_i \in \Gamma$.
In practice, the matrix $\left(\mathbf{I} - \dfrac{\boldsymbol{\Xi} \boldsymbol{\Xi}^\intercal}{\boldsymbol{\Xi}^\intercal \boldsymbol{\Xi}} \right)$ does not need to be explicitly assembled. Noting that $\dfrac{\boldsymbol{\Xi}^\intercal \mathbf{p}}{\boldsymbol{\Xi}^\intercal \boldsymbol{\Xi}}$ actually
calculates the mean of the entities in $\mathbf{p}$, 
the application of $\left(\mathbf{I} - \dfrac{\boldsymbol{\Xi} \boldsymbol{\Xi}^\intercal}{\boldsymbol{\Xi}^\intercal \boldsymbol{\Xi}} \right)$ can be implemented as: first calculating the mean of all entities of the vector $\mathbf{p}$ in parallel and then subtracting the mean from each entity of $\mathbf{p}$. By such, the block structure of $\mathbf{F}$ can be preserved, and the parallel scalability is unaffected.


\section{Numerical Examples}\label{sec:numerical_experiment}
In this section, we systematically assess the effectiveness and scalability of our proposed monolithic GMG preconditioner through several numerical examples, including pure fluid flows and fluid-solid interaction problems.

\subsection{Pure fluid flows}\label{subsec:2d_pure_fluid}
We first verify our GMLS discretization, GMG preconditioner, and parallel implementation by solving problems of pure fluid flows. We start with a simple Taylor–Green vortex flow followed by a more complicated Stokes flow in an artificial vascular network. 

\subsubsection{Taylor–Green vortex}\label{subsec:Taylor-Green}
The fluid domain is set as $\Omega_f = [-1, 1] \times [-1, 1]$. Given the source term in Eq. \eqref{eq:governing_eq_recast} as:
\begin{equation}\label{eq:Taylor_Green_source}
    \mathbf{f} =
    \begin{bmatrix}
        2 \pi^2 \cos \left( \pi x \right) \sin \left( \pi y \right)  + 2\pi \sin (2\pi x) \\
        -2 \pi^2 \sin \left( \pi x \right) \cos \left( \pi y \right) + 2\pi \sin (2\pi y) 
    \end{bmatrix}\;, \quad \forall \ \mathbf{x}=(x,y) \in \Omega_f \;,
\end{equation}
and the no-slip BC for velocity as
\begin{equation}
    \left \{
    \begin{aligned}
        u = & \cos \left( \pi x \right) \sin \left( \pi y \right)  \\
        v = & - \sin \left( \pi x \right) \cos \left( \pi y \right)
    \end{aligned}
    \right. \;,
    \quad \forall \ \mathbf{x} = (x, y) \in \Gamma_0 \;,
    \label{eq:green_vortex_2d_bc}
\end{equation}
where $\Gamma_0 = \partial \Omega_f$ denotes the outer boundary of the fluid domain, the analytical solution of Eq. \eqref{eq:governing_eq_recast} is the following:
\begin{equation}
    \left \{
    \begin{aligned}
        u = & \cos \left( \pi x \right) \sin \left( \pi y \right)   \\
        v = & - \sin \left( \pi x \right) \cos \left( \pi y \right) \\
        p = & -\cos (2 \pi x) - \cos (2 \pi y)
    \end{aligned}
    \right. \;,
    \quad \forall \ \mathbf{x}=(x,y) \in \Omega_f.
    \label{eq:green_vortex_2d}
\end{equation}

\begin{figure}[htp]
    \centering
    \includegraphics[width=8cm]{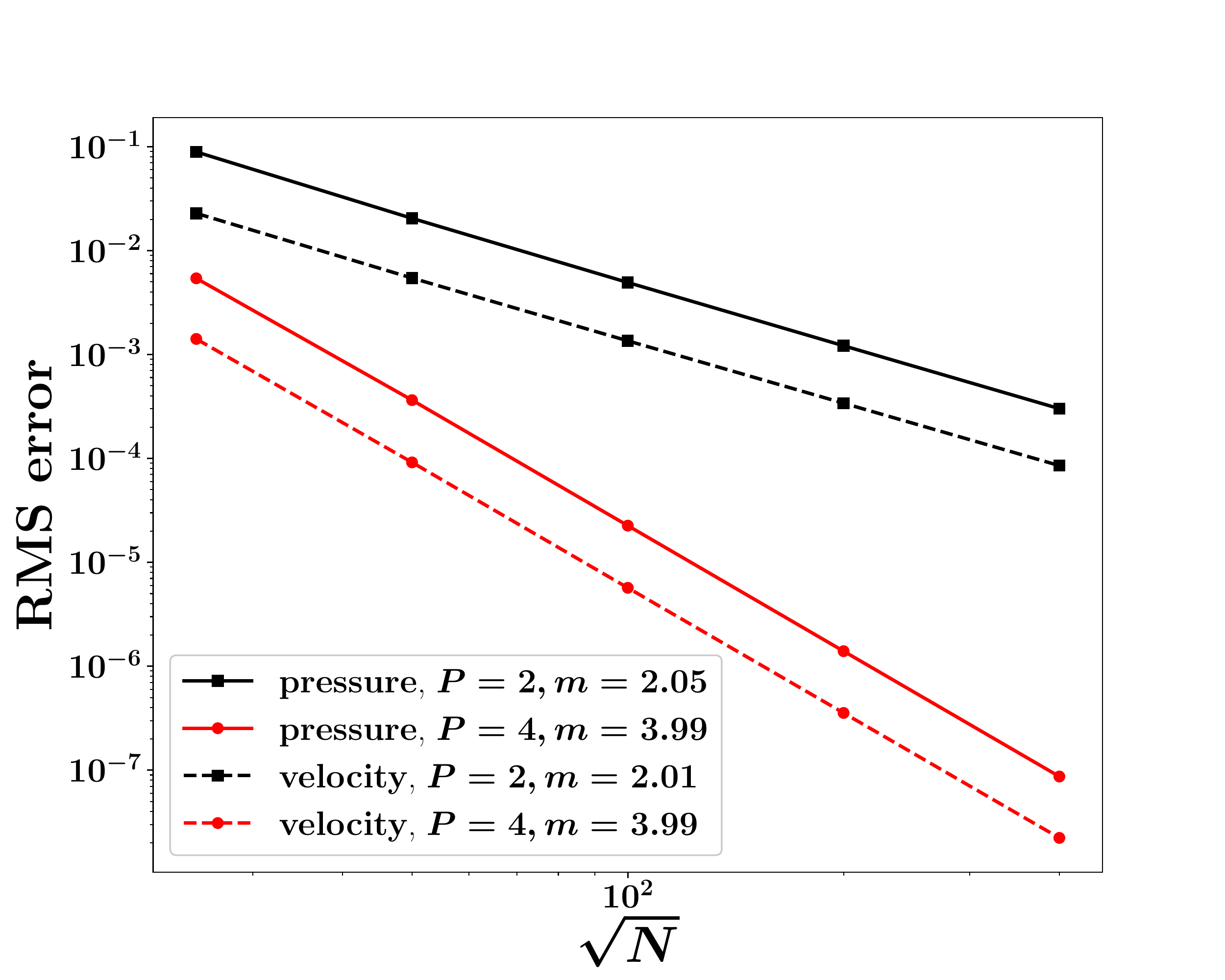}
    \caption{Pure fluid flow$-$Taylor–Green vortex: RMS errors and convergence for the numerical solutions of the velocity (dashed line) and pressure (solid line). Here, $N$ denotes the total number of GMLS nodes; $P$ denotes the order of polynomial basis used in the GMLS discretization; $m$ is the slope of each line.} 
    \label{fig:pure_flow_RMS}.
\end{figure}

Our numerical solutions of this problem obtained with the 2nd-order or 4th-order GMLS discretization, respectively, are compared with the analytical solution. The root mean square (RMS) errors are computed for both velocity and pressure. As shown in \Cref{fig:pure_flow_RMS}, the numerical results exhibit the consistent 2nd-order or 4th-order convergence, as theoretically expected; the velocity and pressure fields achieve equal-order optimal convergence.

After verifying the accuracy and convergence of the numerical solutions, we next examine the scalability of the proposed GMG preconditioner, as well as the parallel scalability of our implementation. To this end, the numbers of CPU cores $N_c$ are varied from $25$ to $289$ in the tests. For such a pure fluid flow problem without singularities, adaptive refinement is not needed, and hence only uniform refinement is conducted. To maintain the same number of GMLS nodes distributed in each CPU core, we start from $1\times1$ GMLS nodes in each core, and then execute the same uniform refinement for the GMLS nodes within each core. After seven iterations of uniform refinement, there are $128 \times 128$ GMLS nodes in each CPU core. 
\begin{figure}[htp]
    \centering
    \begin{subfigure}{.47\textwidth}
        \centering
        \includegraphics[width=.98\textwidth]{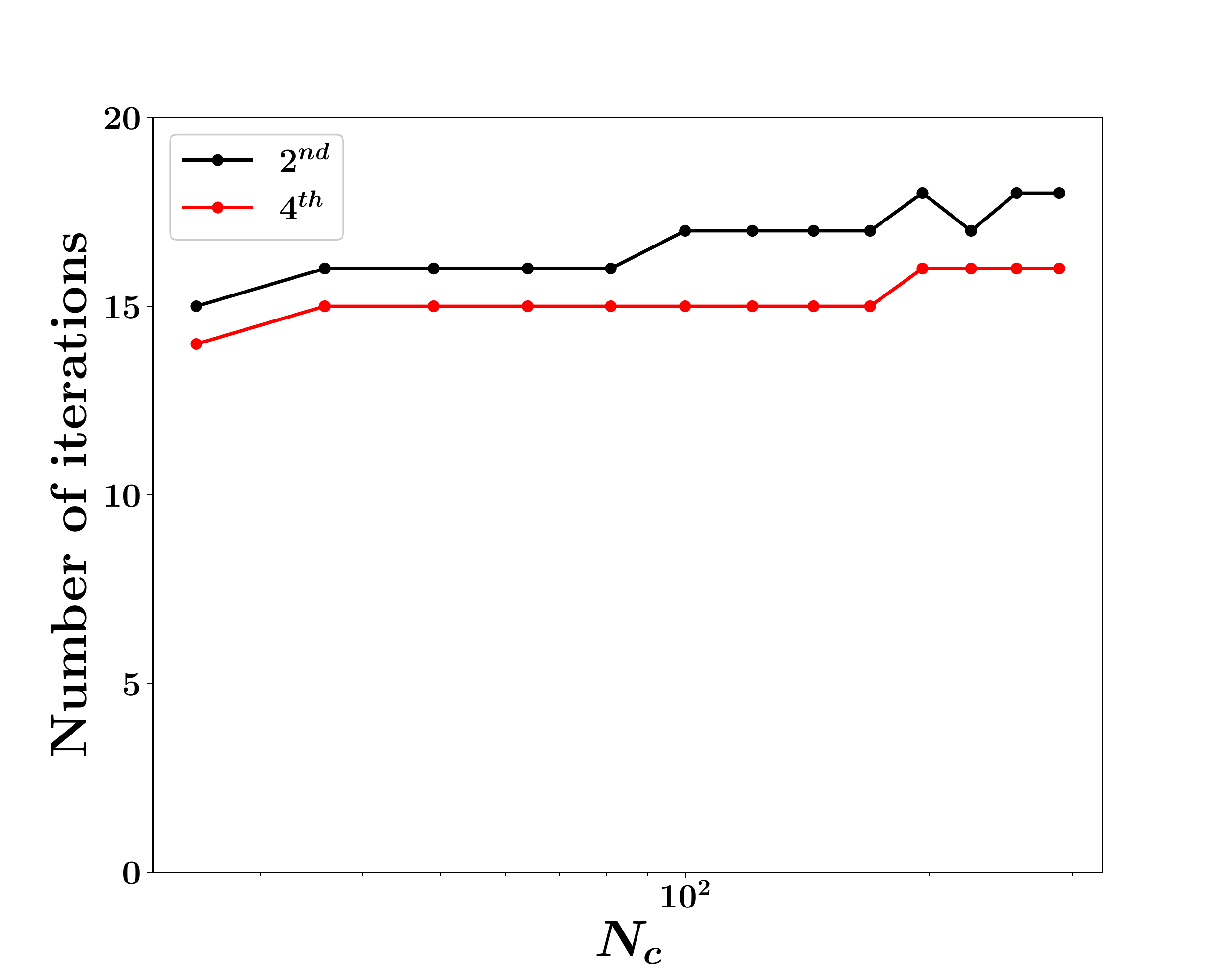}
        \caption{Number of iterations required for the linear solver to converge.}
        \label{fig:pure_flow_ite_num}
    \end{subfigure}
    \quad \quad
    \begin{subfigure}{.47\textwidth}
        \centering
        \includegraphics[width=.98\textwidth]{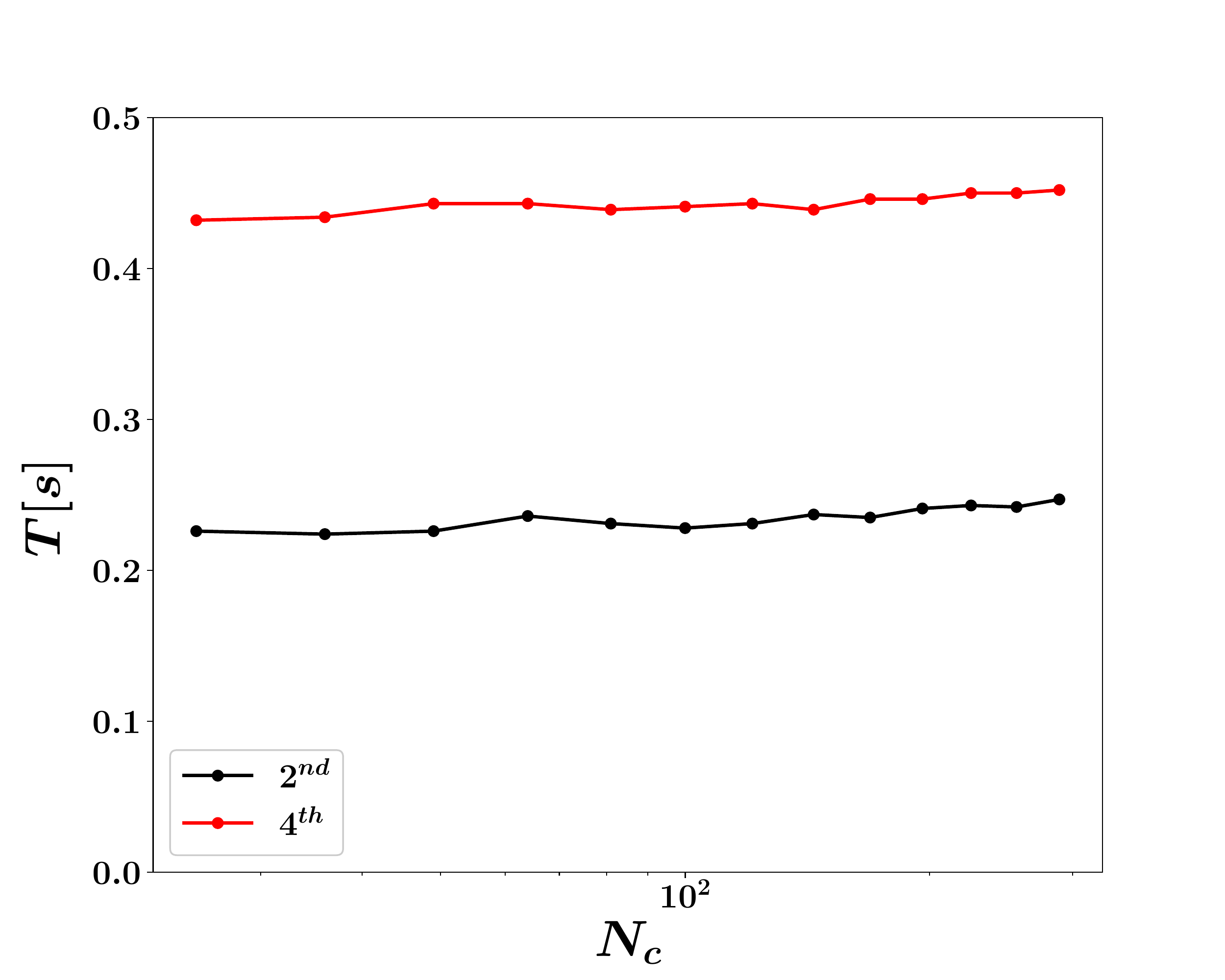}
        \caption{Computer time (in seconds) spent for a single step of preconditioning.}
        \label{fig:pure_flow_pcapply}
    \end{subfigure}
    \caption{Pure fluid flow$-$Taylor–Green vortex: Scalability of the proposed GMG preconditioner and the weak scalability of our parallel implementation of the preconditioner, tested for different order of GMLS discretization (black for the 2nd order and red for the 4th order). Here, $N_c$ denotes the number of CPU cores used in each test.}\label{fig:pure_flow_linear_system}
\end{figure}
To verify the scalability of the proposed GMG preconditioning method, we examine the number of iterations required for the linear solver to converge. To evaluate the parallel scalability of our implementation, we record the computer time spent for a single step of preconditioning. Our tests are particularly for the linear system generated from the last iteration ($I = 7$) of (uniform) refinement, i.e., with the most DOFs. In \Cref{fig:pure_flow_ite_num}, we show the variation of the number of iterations with respect to different numbers of CPU cores used and find that it only varies slightly from $15$ to $18$ for the $2$-nd order GMLS discretization and $14$ to $16$ for the $4$-th order GMLS as the number of CPU cores increases from $25$ to $289$. The computer time spent on a single step of preconditioning stays almost constant with increasing CPU cores, independent of the order of GMLS discretization, as depicted in \Cref{fig:pure_flow_pcapply}. By these results, we demonstrate the scalability of the proposed GMG preconditioner and the weak scalability of our parallel implementation of the preconditioner.

\subsubsection{Artificial vascular network}
\begin{figure}[htp]
    \centering
    \includegraphics[height=8cm]{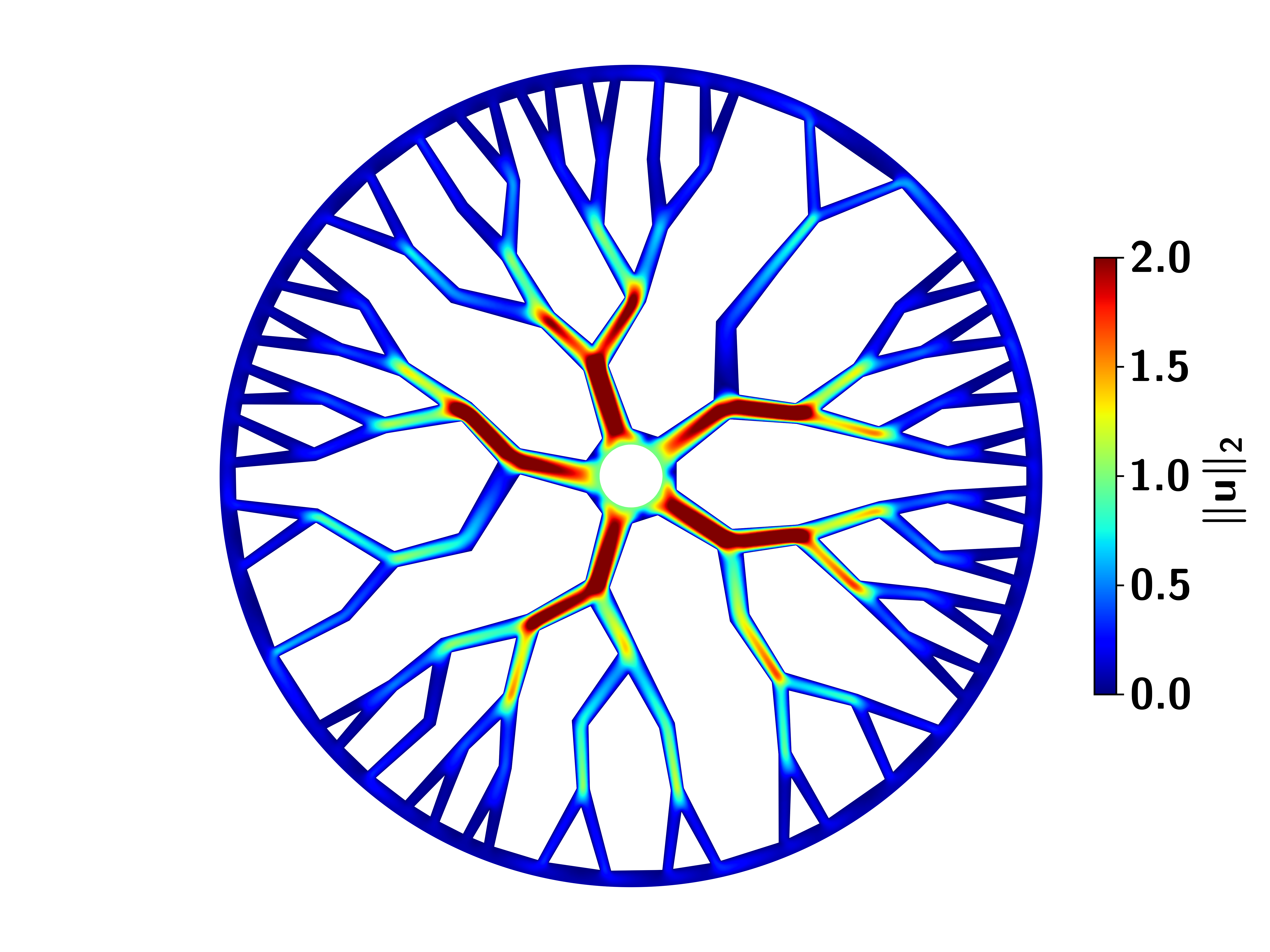}
    \caption{Pure fluid flow$-$Artificial vascular network: Computed velocity field, where the color bar indicates the magnitude of velocity.}
    \label{fig:artificial_vascular_network}
\end{figure}

We next solve a more complicated pure fluid flow problem, which is in an artificial vascular network, mimicking the structure of a zebrafish's eye \cite{alvarez2007genetic}. The artificial vascular network is built as in \Cref{fig:artificial_vascular_network}. A constant inflow flux drives the flow in this network at the inner circular boundary in the center, and a constant outflow flux is imposed at the outermost circular boundary. The overall volume of the fluid in the network is conserved. All the other boundaries in the network are imposed no-slip BCs for the velocity.  
Due to the irregular computational domain with corner singularities, adaptive refinement is required such that the numerical solution can achieve optimal convergence. Starting from a coarse resolution at $\Delta x^0 = 0.02$ and setting $\alpha = 0.8$ (marking percentage) in \textbf{MARK} stage, we perform $8$ iterations of adaptive refinements. The velocity field computed after the 8th iteration of adaptive refinement is shown in \Cref{fig:artificial_vascular_network}, where there are in total 930,240 GMLS nodes, and the 2nd-order GMLS discretization is employed. 

\begin{figure}[htp]
    \centering
        \includegraphics[width=8cm]{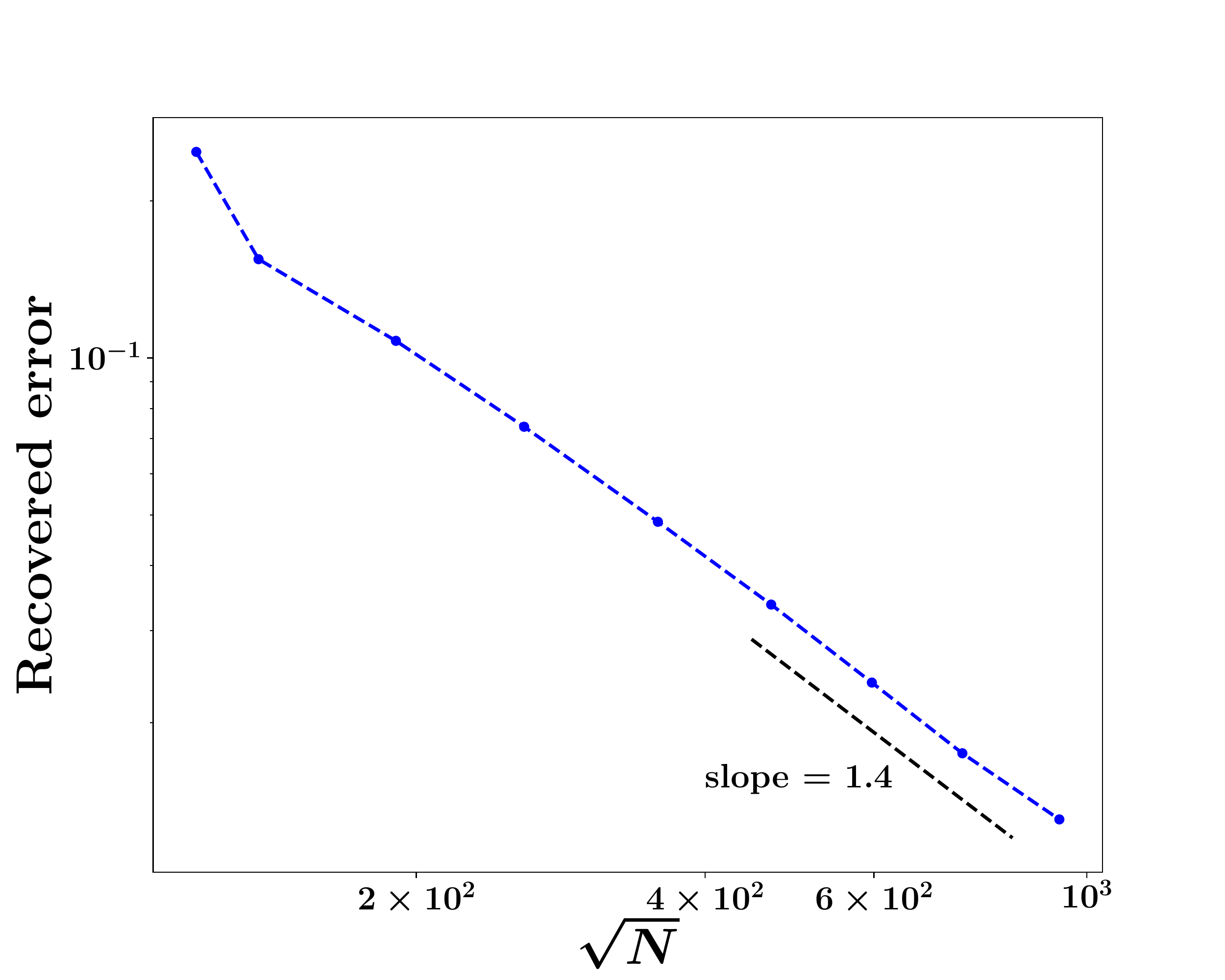}
        \caption{Pure fluid flow$-$Artificial vascular network: Convergence of the total recovered error. Here, the slope is regressed from the last 4 points; $N$ denotes the total number of GMLS nodes.}
        \label{fig:convergence_network}
\end{figure}

\begin{figure}[htp]
    \centering
    \begin{subfigure}{.47\textwidth}
        \centering
        \includegraphics[height=7cm]{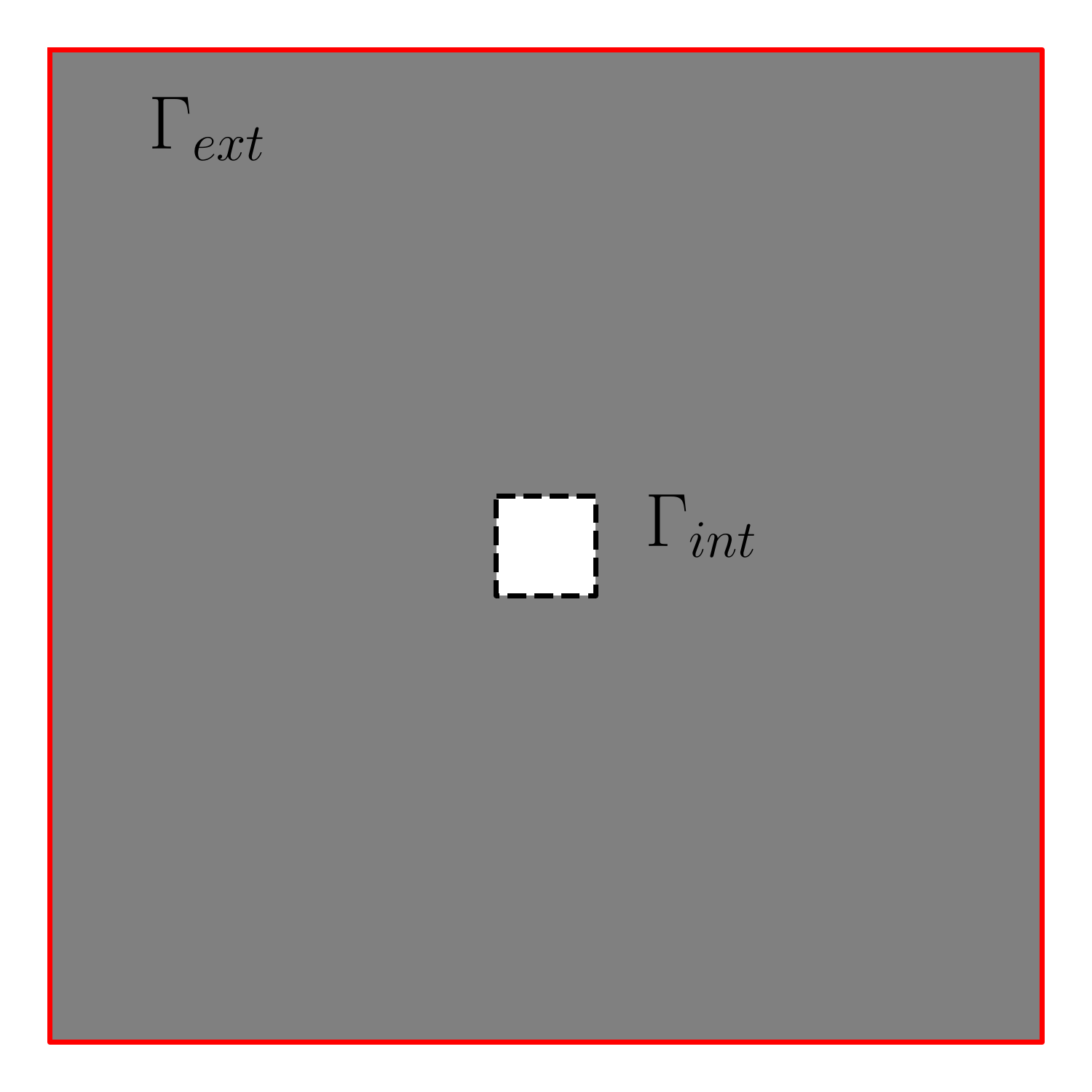}
        \caption{Illustration of the computational domain: Fluid (gray) is confined in  $\Omega_f$ by $\Gamma_{ext}$ and $\Gamma_{int}$. $\Gamma_{ext}$ (red solid line) is a square with side length of 1. $\Gamma_{int}$ (black dashed line) is at the center of domain and varied from a square to a hexagon, triangle, or parallelogram, all with an equal side length of $0.2$.}
        \label{fig:boundary_illustration}
    \end{subfigure}
    \quad \quad
    \begin{subfigure}{.47\textwidth}
        \centering
        \includegraphics[height=7cm]{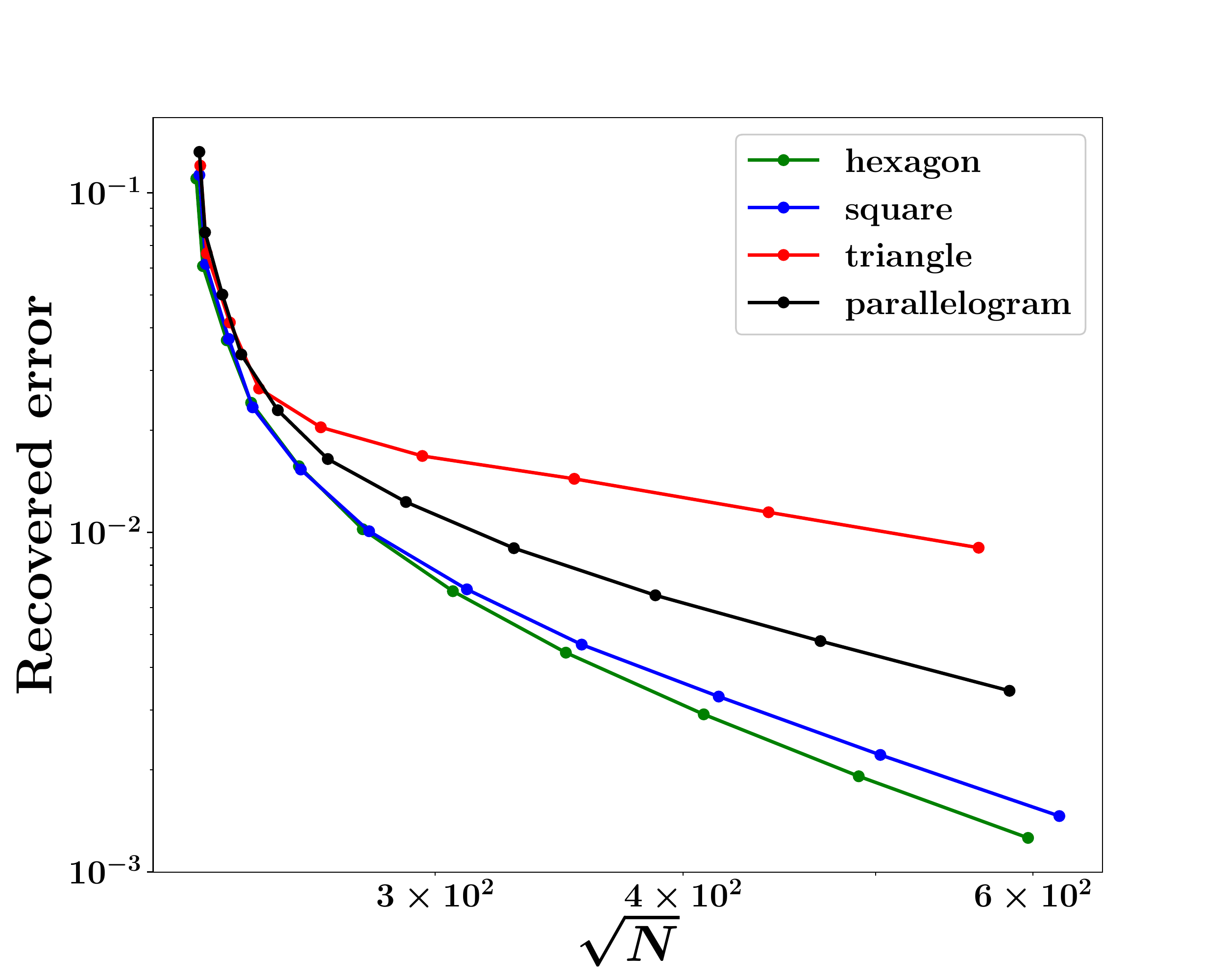}
        \caption{Convergence of the total recovered error. $N$ denotes the total number of GMLS nodes. The convergence rate is estimated as the slope regressed from the last 4 data points: slop$\approx$2.0 for hexagon and square; slop$\approx$1.6 for parallelogram; slop$\approx$1.0 for triangle.}
        \label{fig:error_shape}
        \end{subfigure}
    \caption{Convergence test with respect to the nonconvexity of computational domain.}
\end{figure}

Due to a lack of the true solution, we calculate the total recovered errors (Eq. \eqref{eq:global_posteriori_err}), instead of true errors, during iterations of adaptive refinement to evaluate the accuracy and convergence of the numerical solutions, as shown in \Cref{fig:convergence_network}. 
We note that the convergence rate does not reach the theoretical 2nd order. This arises from the presence of significant portions of nonconvex boundaries in the computational domain. It is well-known that nonconvexity can lead to low regularity of the solutions to elliptic PDEs, which in turn affects the performance and liability of recovery-based \textit{a posteriori} error estimator, resulting in deteriorated convergence in adaptive refinement~\cite{xu2004analysis,susanne2008mathematical}. Therefore, we design a numerical test to elucidate this issue. In the test, we let the computational domain be constrained by two boundaries, an exterior boundary $\Gamma_{ext}$ and an interior boundary $\Gamma_{int}$. $\Gamma_{ext}$ is set the same as that in \S\ref{subsec:Taylor-Green} with the same BCs (Eq. \eqref{eq:green_vortex_2d_bc}). $\Gamma_{int}$ is stationary but with different shapes, including square, hexagon, triangle, and parallelogram. Taking $\Gamma_{int}$ of a square as an example, the computational domain $\Omega_f$ is depicted in \Cref{fig:boundary_illustration}. For each shape of $\Gamma_{int}$, the problem is solved with $10$ iterations of adaptive refinement, and the total recovered errors calculated during adaptive refinement is summarized in \Cref{fig:error_shape}. As expected, the convergence rate decreases from $2$ to $1$ when $\Omega_{int}$ changes from hexagon to triangle, with increasing nonconvexity. This confirms our sub-optimal convergence results observed in \Cref{fig:convergence_network}.

\begin{figure}[htp]
    \centering
    \includegraphics[width=8cm]{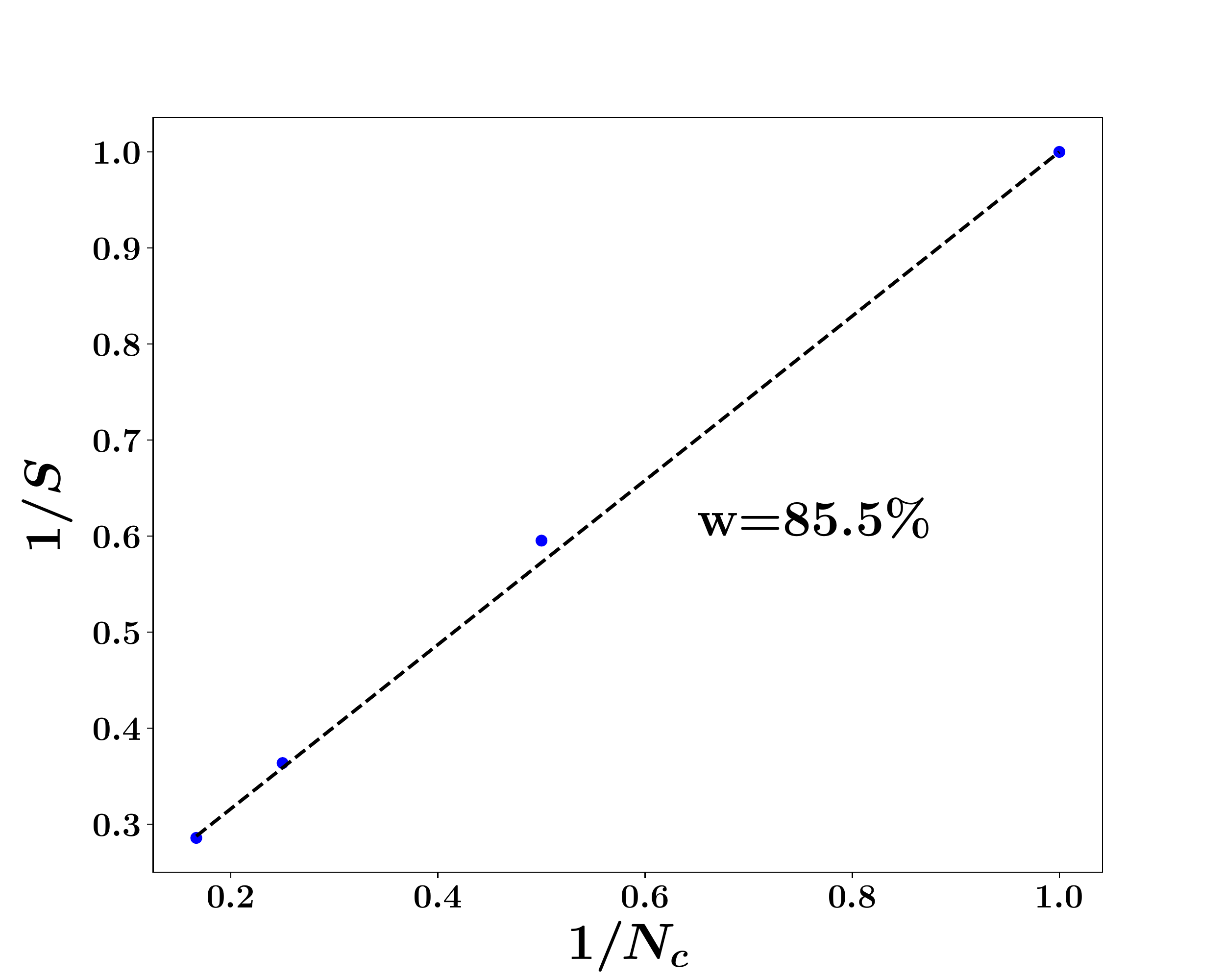}
    \caption{Pure fluid flow$-$Artificial vascular network: Parallel portion $w$ determined from Amdahl's law in Eq. \eqref{eq:reciprocal_Amdahl}.}
    \label{fig:parallel_portion}
\end{figure}

After examining the convergence of our numerical solutions, we further assess the scalability of our proposed GMG preconditioner and parallel implementation. 
Due to the complexity of the computational domain, we cannot guarantee that the total number of GMLS nodes resulted in each adaptive refinement iteration increases proportionally, 
thus, the weak scaling test is not appropriate for this problem. Instead, we perform a strong scaling test to demonstrate the parallel scalability of the proposed GMG preconditioning method. In the strong scaling test, the number of cores invoked changes from 40 to 240 cores in the simulation. The statistics of the test is collected in \Cref{tab:strong_scaling_artificial_vascular}.  The average DOFs per core after the $8$-th adaptive refinement iteration is listed in the column \emph{DOFs/Core}, which shows the scale of the workload to this problem. As expected, the workload decreases when the number of cores increases, since we are doing a strong scalability test and the total DOFs are fixed. The average number of GMRES iterations required in each adaptive refinement iteration step is listed in the column \emph{Average Iterations}. It stays around 66, independent of different numbers of cores invoked, which indicates that our parallel implementation does not affect the convergence of the proposed GMG preconditioning method. To evaluate the parallel scalability of our implementation, we record the CPU wall time spent for solving the linear systems in each adaptive refinement iteration and sum them up, which is denoted as $T_s$ in \Cref{tab:strong_scaling_artificial_vascular}; we also track the total CPU wall time spent for the entire simulation, denoted as $T$. 
\begin{table}[ht]
    \centering
    \caption{Pure fluid flow$-$Artificial vascular network: Strong scaling test.}
    \label{tab:strong_scaling_artificial_vascular}
    \begin{tabular}[t]{rcccccc}
        \toprule
        $N_c$ & DOFs/Core & Average Iterations & Time for \eqref{eq:structure_linear_system} $T_s$ [s] & Overall Time $T$ [s] & Overall Speedup $S$ \\
        \midrule
        40 & 69,768 & 65.5 & 90.52 & 128.47 & 1.00 \\
        80 & 34,884 & 67.1 & 52.43 & 76.50 & 1.68 \\
        160 & 17,442 & 66.5 & 33.71 & 46.73 & 2.75 \\
        240 & 11,628 & 66.3 & 28.29 & 36.68 & 3.50 \\
        \bottomrule
    \end{tabular}
\end{table}
By comparing $T_s$ and $T$, we note that the time spent for solving the linear systems indeed dominates the overall computing time, by more than $70\%$,
regardless of how many cores used. The GMLS discretization can be trivially parallelized, and hence the overall parallel scalability is dictated by the proposed preconditioning method. By invoking different numbers of cores, we can compare the overall execution time and thereby evaluate the speed up factor (denoted as $S$) using the least number of cores (i.e., 40 cores) as the base. 
We hence use Amdahl's law \cite{bryant2003computer} to assess the performance of our parallel implementation, which is given by:  
\begin{equation*}
    S = \dfrac{1}{(1 - w) + \dfrac{w}{N_c}}\;,
\end{equation*}
where $w$ denotes the parallel portion of a numerical solver. To estimate the value of $w$ for our solver, we take the reciprocal of Amdahl's law as: 
\begin{equation}\label{eq:reciprocal_Amdahl}
    \dfrac1S = (1 - w) + \dfrac{w}{N_c}\;.
\end{equation}
Using the data of $S$ in \Cref{tab:strong_scaling_artificial_vascular}, we can determine $w$ by linear regression. In particular, we employ only the data of $N_c=40$ and $N_c=80$ to determine $w$, and then use the data of $N_c=160$ and $N_c=240$ for testing. As depicted in \Cref{fig:parallel_portion}, the last two data points fall very close to the line fitted from the first two data points, and the estimated parallel portion $w = 85.5\%$. We hence demonstrate the consistency and scalability of our parallel implementation of the proposed monolithic GMG preconditioner. 

\subsection{Fluid-solid interactions}
After pure fluid flows, we next address fluid-solid interaction problems with the inclusion of multiple solid bodies of different shapes. 

\subsubsection{Duplicate cells with cylinders} 
\begin{figure}[htp]
    \centering
    \includegraphics[width=8cm]{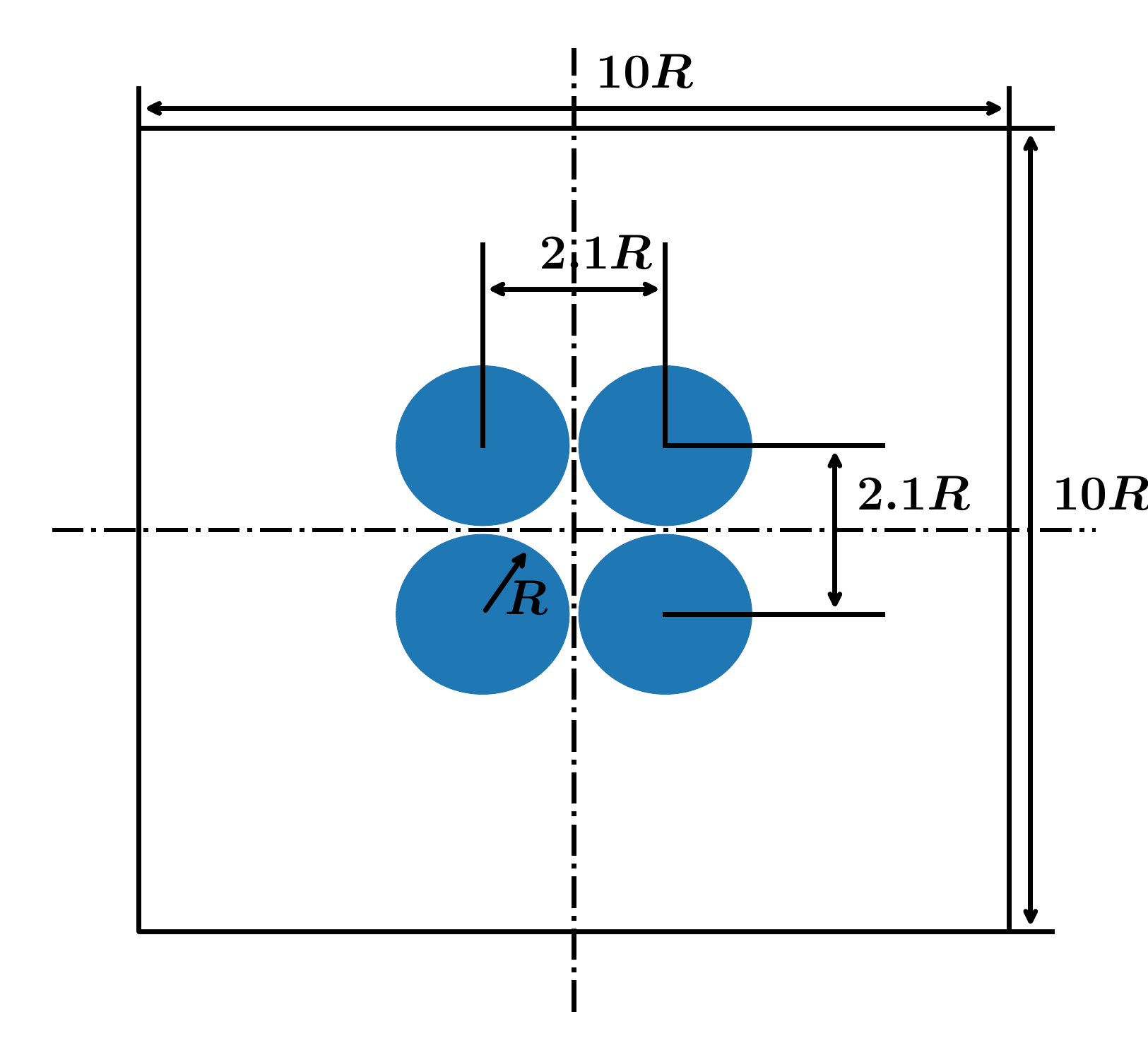}
    \caption{Fluid-solid interactions$-$Duplicate cells with cylinders: Schematic of a square cell with four cylinders.}
    \label{fig:multi_solid_scaling_single_cell_illustration}
\end{figure}

First, for an accurate assessment of the scalability of our proposed monolithic GMG preconditioner, we design a numerical example that allows both the solid bodies and the GMLS nodes to be evenly distributed among computer cores. To this end, the entire computational domain is partitioned into $N_c$ square cells, each of which includes four cylinders, as shown in \Cref{fig:multi_solid_scaling_single_cell_illustration}. 
All DOFs associated in each square cell are allocated to a single CPU core. Thus, there are in total $N_s=4N_c$ solid bodies, where $N_c$ is the number of CPU cores for each test. We intentionally place the cylinders in close contacts to demonstrate that our spatially adaptive GMLS method can resolve the singularities governing the lubrication effects. To maintain the same total DOFs after adaptive refinement in each cell, the four cylinders are placed near the center of each cell such that the distances between the solid bodies in different cells are much larger than the distances between the solid bodies within one cell.   

\begin{figure}[htp]
    \centering
    \includegraphics[width=10cm]{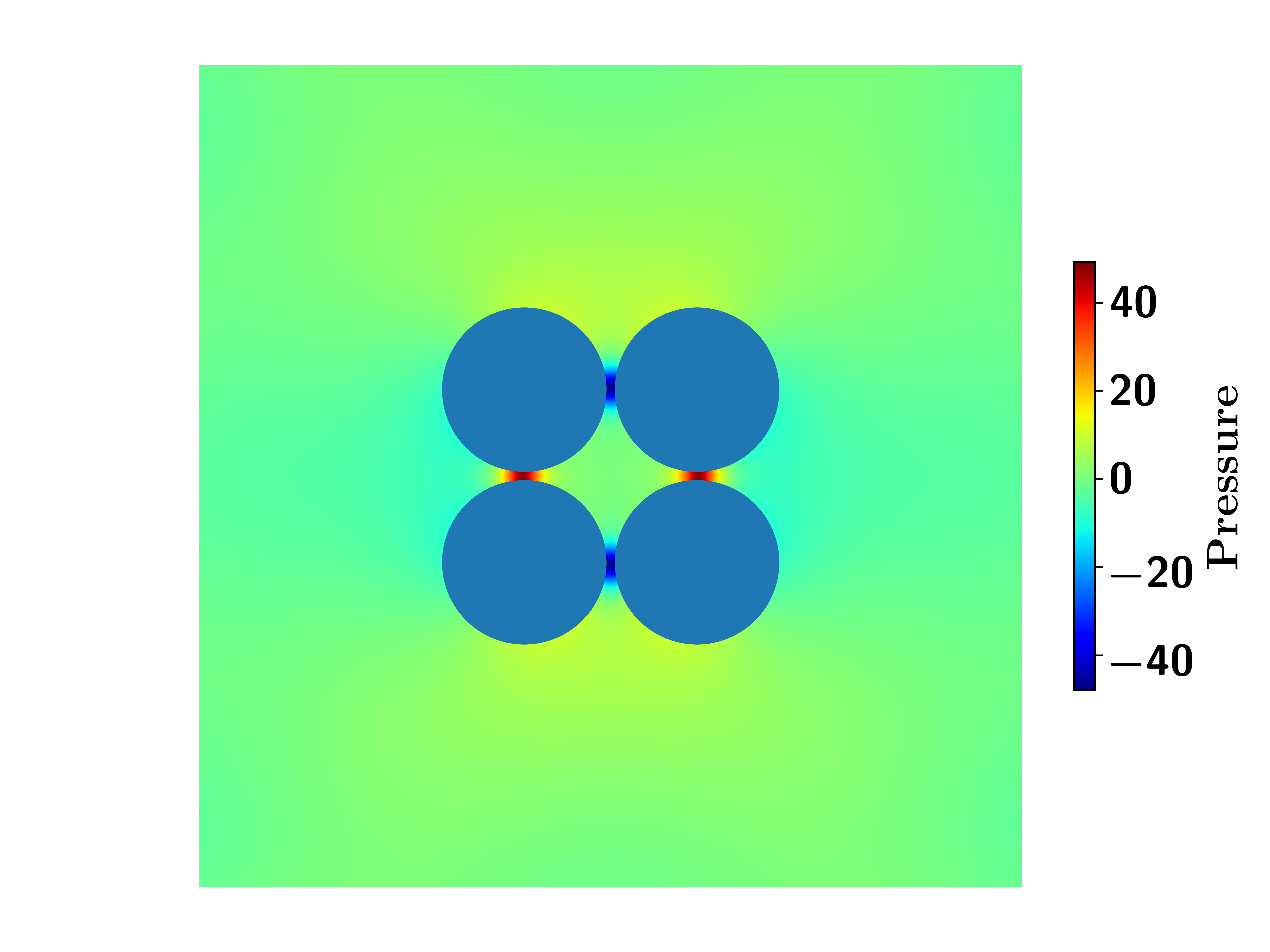}
    \caption{Fluid-solid interactions$-$Duplicate cells with cylinders: The pressure field computed in each cell. The color is correlated to the magnitude of pressure.}
    \label{fig:single_cell_pressure}
\end{figure}

The source term in Eq. \eqref{eq:governing_eq_recast} and the BC at the outer boundary of the entire computational domain are set the same as those in \S\ref{subsec:Taylor-Green}, i.e., Eqs. \eqref{eq:Taylor_Green_source}- \eqref{eq:green_vortex_2d_bc}. The 2nd-order GMLS discretization is employed for solving this problem. Seven iterations of adaptive refinement with $\alpha = 0.8$ are conducted until the total recovered error reaches the preset error tolerance $e_{tol} = 10^{-3}$ in Eq. \eqref{eq:error_tolerance}. The stopping criterion for the GMRES iteration is set as $10^{-6}$. The resultant pressure field in a single cell is shown in \Cref{fig:single_cell_pressure}. We can see that all the singular pressures within the narrow gaps between cylinders are correctly captured. 

\begin{figure}[htp]
    \centering
    \begin{subfigure}{.47\textwidth}
        \centering
        \includegraphics[width=.98\textwidth]{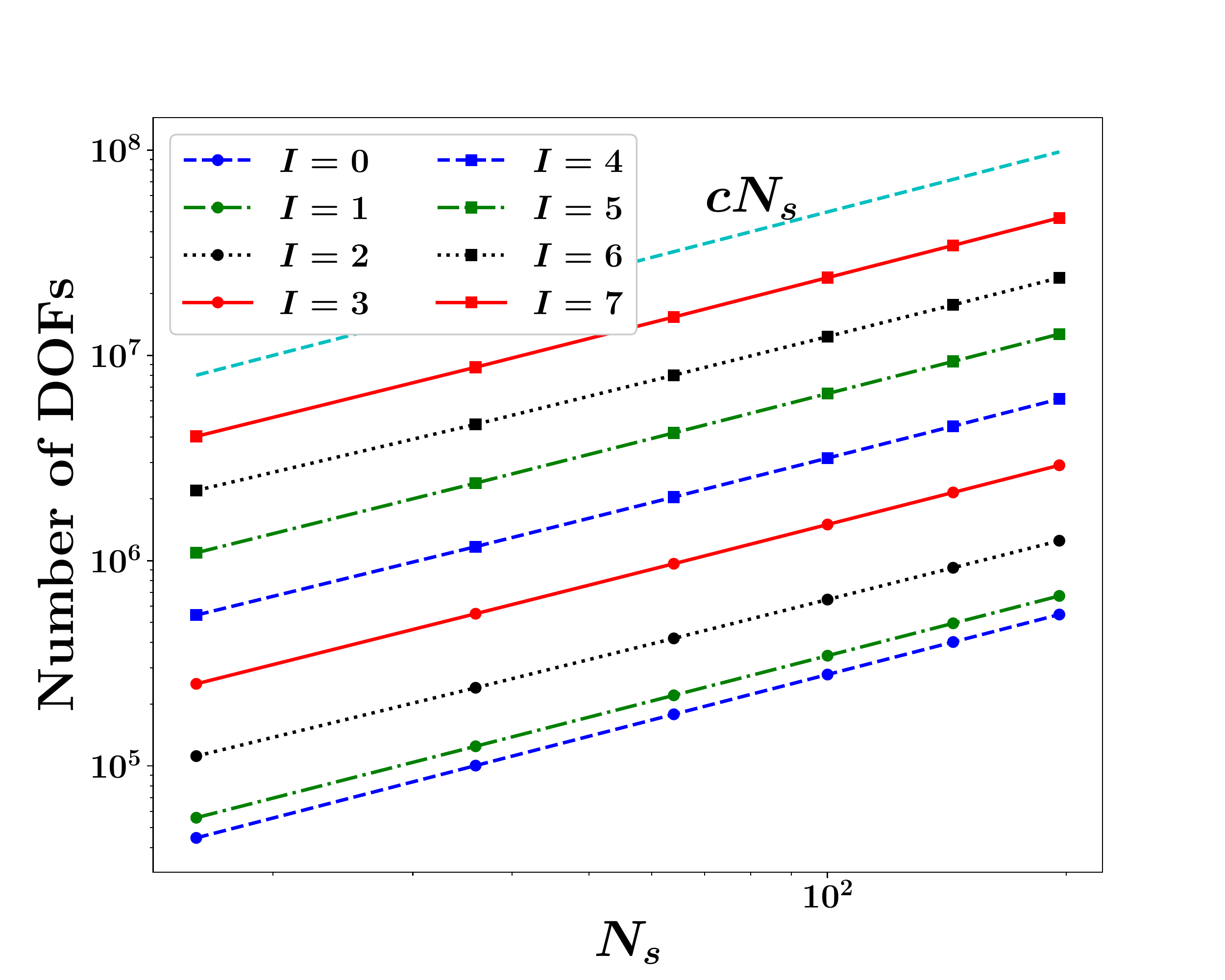}
        \caption{Number of total DOFs.}
        \label{fig:multi_solid_dof_refinement}
    \end{subfigure}
    \quad \quad
    \begin{subfigure}{.47\textwidth}
        \centering
        \includegraphics[width=.98\textwidth]{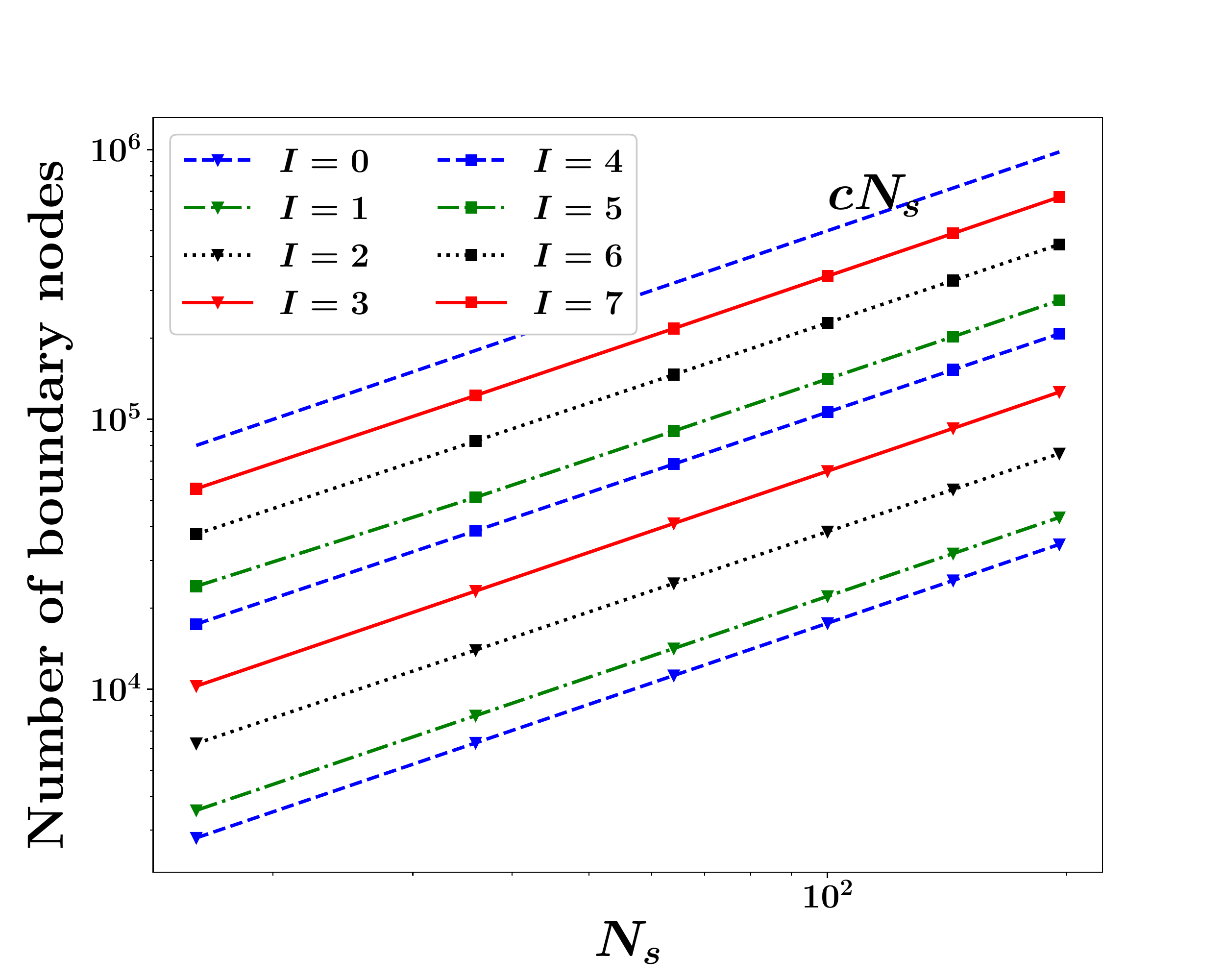}
        \caption{Number of the boundary GMLS nodes.}
        \label{fig:surface_particle_num}
    \end{subfigure}
    \caption{Fluid-solid interactions$-$Duplicate cells with cylinders: The growths of total DOFs and the boundary (GMLS) nodes with respect to 
    the number of solids included in the domain and different iterations of adaptive refinement.}
\end{figure}

Before we assess the parallel scalability of our preconditioner, we examine how the total DOFs consisting of fluid (GMLS) nodes and boundary (GMLS) nodes grow with adaptive refinement and the inclusion of more solid bodies. In \Cref{fig:multi_solid_dof_refinement}, we can see that the total DOFs increases linearly with respect to the number of solids included in the domain, regardless of after any iteration of adaptive refinement; in \Cref{fig:surface_particle_num}, we find that the total boundary nodes also increases linearly with respect to the number of solids, after each iteration of adaptive refinement. From these examinations, we confirm that the workload related to the node-wise GS smoother $\mathbf{S}_{\mathbf{F}}$, the additive Schwarz-type smoother $\mathbf{S}_{\mathbf{N}}$, and the MG preconditioner can be evenly distributed among all cores. Only with that ensured, we can make an accurate assessment of scalability.

\begin{figure}[htp]
	\centering
	\includegraphics[width=8cm]{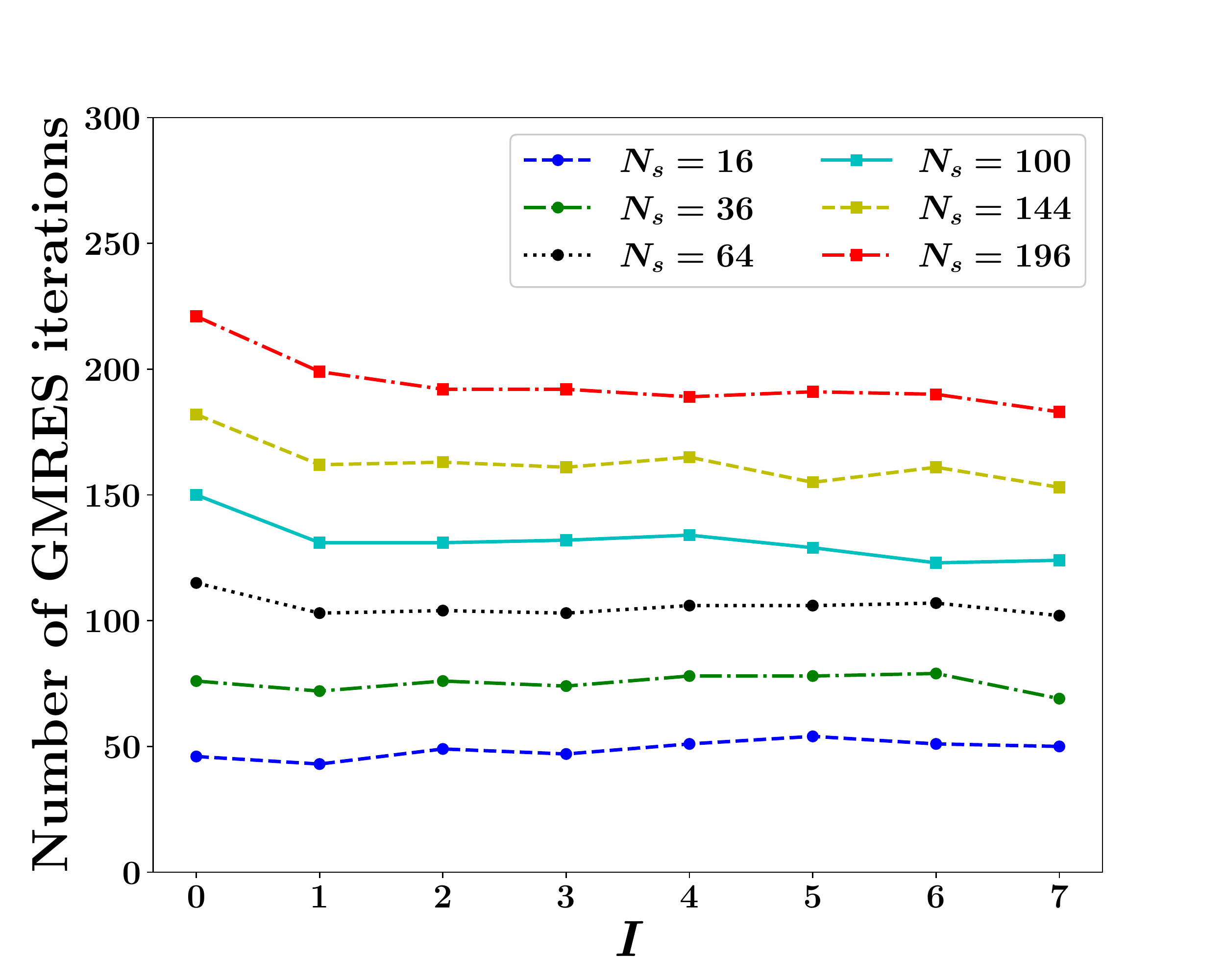}
	\caption{Fluid-solid interactions$-$Duplicate cells with cylinders: The number of GMRES iterations required in the \textbf{SOLVE} stage of each adaptive refinement iteration, for fixed numbers of solids.}\label{fig:multi_solid_ite_refinement}
\end{figure}

\begin{figure}[htp]
    \centering
    \begin{subfigure}{.47\textwidth}
        \centering
        \includegraphics[height=7cm]{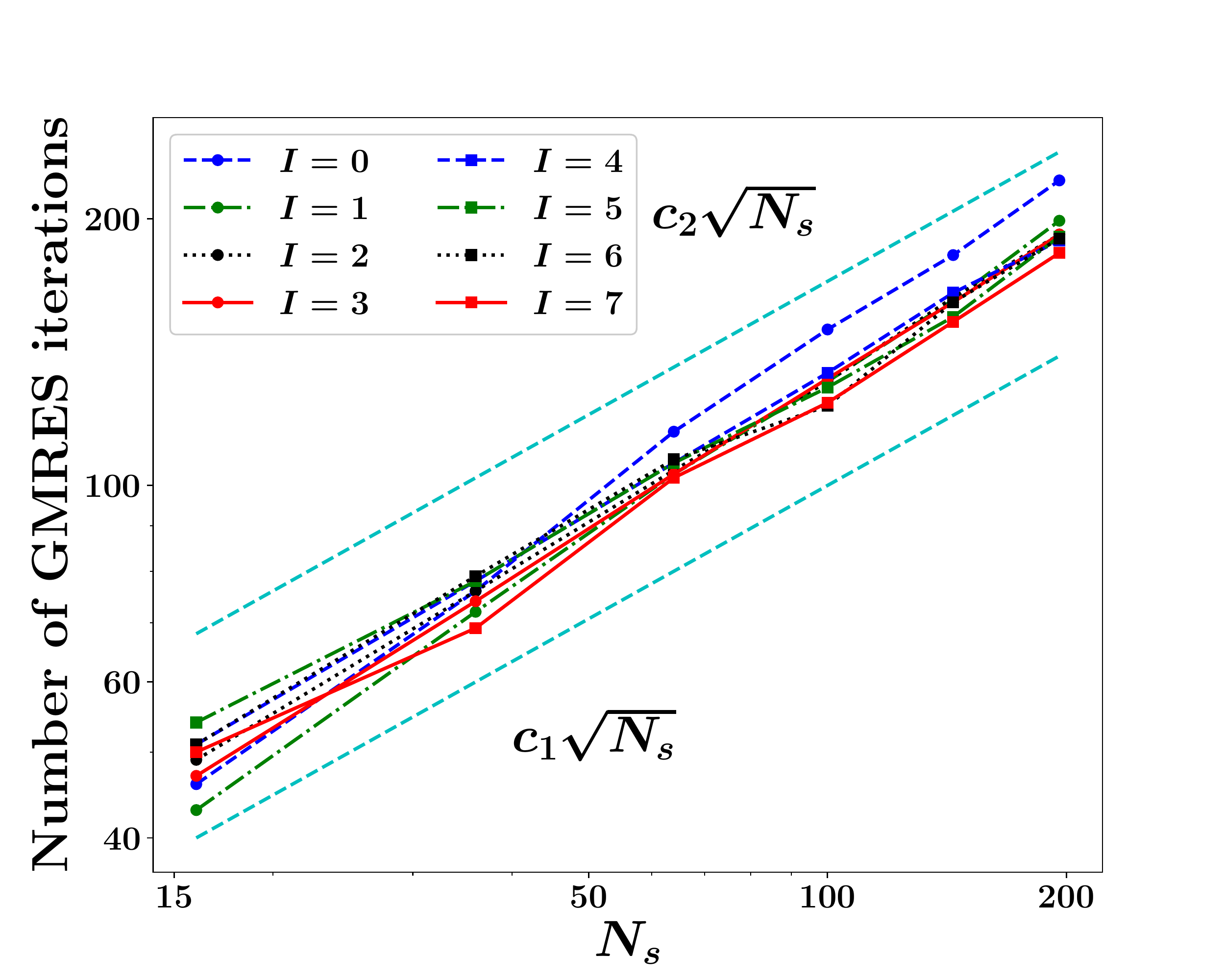}
        \caption{Scaling of the number of GMRES iterations required at different iteration step of adaptive refinement.}\label{fig:multi_solid_ite_Nc}
    \end{subfigure}
    \quad \quad
    \begin{subfigure}{.47\textwidth}
        \centering
        \includegraphics[height=7cm]{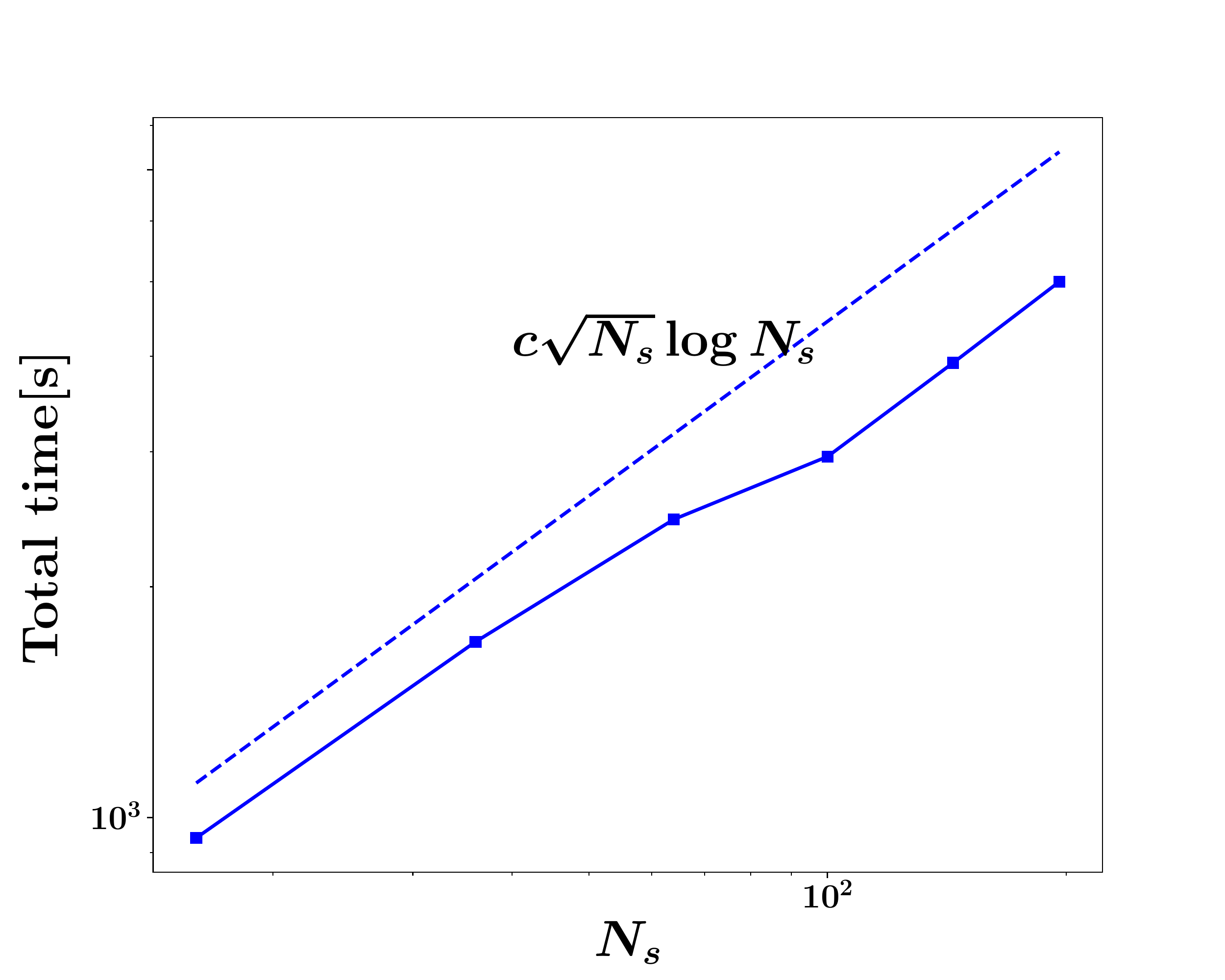}
        \caption{Scaling of the total computer time spent for all seven adaptive refinement iterations.\quad \quad}
        \label{fig:total_time_solids}
    \end{subfigure}
    \caption{Fluid-solid interactions$-$Duplicate cells with cylinders: Scalability results. Here, $N_s$ denotes the number of solids, and $N_s=4N_c$ with $N_c$ the number of CPU cores.}
    \label{fig:Duplicate_cells_scalability}
\end{figure}

We first check the number of iterations required for the GMRES iterative solver to converge, i.e., to reach the stopping criterion. In~\Cref{fig:multi_solid_ite_refinement}, we show the number of GMRES iterations required in the \textbf{SOLVE} stage of each iteration of adaptive refinement. 
As can be seen, for a fixed number of solids, the number of GMRES iterations required generally stays constant or decreases slightly. Noting that the total DOFs continuously increase during iterations of adaptive refinement, we hence demonstrate the weak scalability of our preconditioner with respect to the number of DOFs for fixed number of solids. With inclusion of more solids, we expect that the number of GMRES iterations required would increase. However, the scaling of this increase is critical for the sake of scalability. \Cref{fig:multi_solid_ite_Nc} shows that the number of GMRES iterations scales with $O(\sqrt{N_s})$, and the scaling is generally consistent for different iteration step of adaptive refinement.
We next check the computer time spent finishing all seven adaptive refinement iterations. Although it includes the time spent on the GMLS discretization and the time for solving the linear system and executing other stages of adaptive refinement, solving the linear system dominates the computer time. \Cref{fig:total_time_solids} depicts how the computer time varies with an increasing number of solids. Overall, it exhibits a scaling of $\mathcal{O}(\sqrt{N_s} \log N_s)$, for which the factor $\sqrt{N_s}$ arises from the scaling of the number of GMRES iterations and the factor $\log N_s$ is mainly contributed by the additional operation introduced in \S\ref{subsec:data_storage} such as Eq.~\eqref{eq:rb_mat}. Note that in this numerical example, the number of CPU cores $N_c$ is proportional to the number of solids $N_s$, in fact $N_s=4N_c$. Thus, the above scaling behaviors are shown in \Cref{fig:Duplicate_cells_scalability} also hold with respect to the number of cores $N_c$.

\subsubsection{Particulate suspensions}
Finally, we perform long-time simulations to examine the performance of the proposed monolithic GMG preconditioning method. In particular, we simulate suspension flows of freely moving particles of different shapes. The flow is driven by the source term in Eq. \eqref{eq:governing_eq_recast} and the BC at the outer boundary of fluid domain ($[-1, 1] \times [-1, 1]$) as in Eqs. \eqref{eq:Taylor_Green_source}- \eqref{eq:green_vortex_2d_bc}. One hundred solid particles are suspended in the flow, subject to bidirectional hydrodynamic couplings. Initially, all particles are evenly distributed throughout the domain. Due to hydrodynamic couplings, particles are freely moving with the flow. The physical time of the entire simulation is $T = 5$. The 2nd-order GMLS is employed for spatial discretization. The 5th-order Runge-Kutta integrator with adaptive time stepping is used for temporal integration to update the particles' translational and angular positions. An initial time step $\Delta t = 0.2$ is applied. In each time step, adaptive $h$-refinement with $\alpha = 0.8$ are conducted until the total recovered error reaches the preset error tolerance $e_{tol} = 10^{-3}$ in Eq. \eqref{eq:error_tolerance}. For the linear solver, the stopping criterion for the GMRES iteration is set as $10^{-6}$.

\subsubsection*{100 similar particles} 
\begin{figure}[htp]
    \centering
    \includegraphics[width=.98\textwidth]{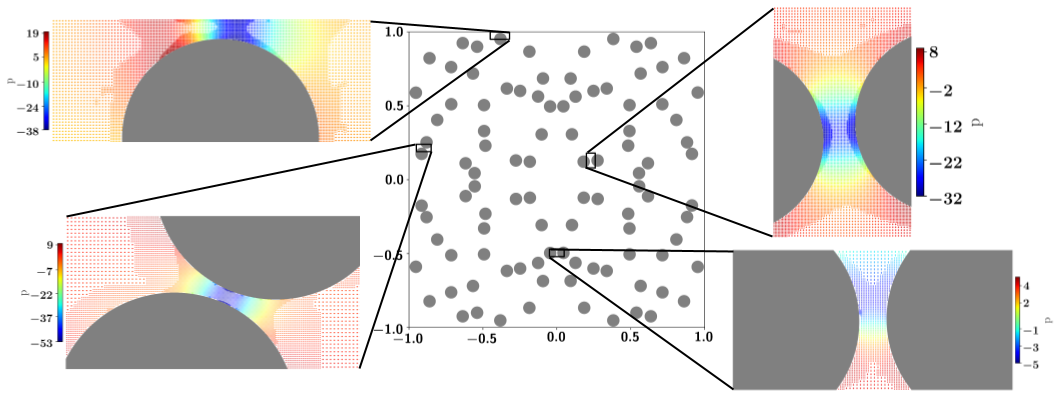}
    \caption{Fluid-solid interactions$-$Particulate suspensions: Configuration of 100 freely moving circular particles in a Taylor–Green vortex flow at the terminal time. The zoom-in images are the computed pressure distributions at selected locations when the particles are either in close contact with each other or with the outer wall, where the color is correlated to the magnitude of pressure, and the point clouds are the GMLS nodes with adaptive refinement.}
    \label{fig:initial_configuration_100_circular}
\end{figure}

\begin{figure}[htp]
    \centering
    \begin{subfigure}{.47\textwidth}
        \includegraphics[height=6cm]{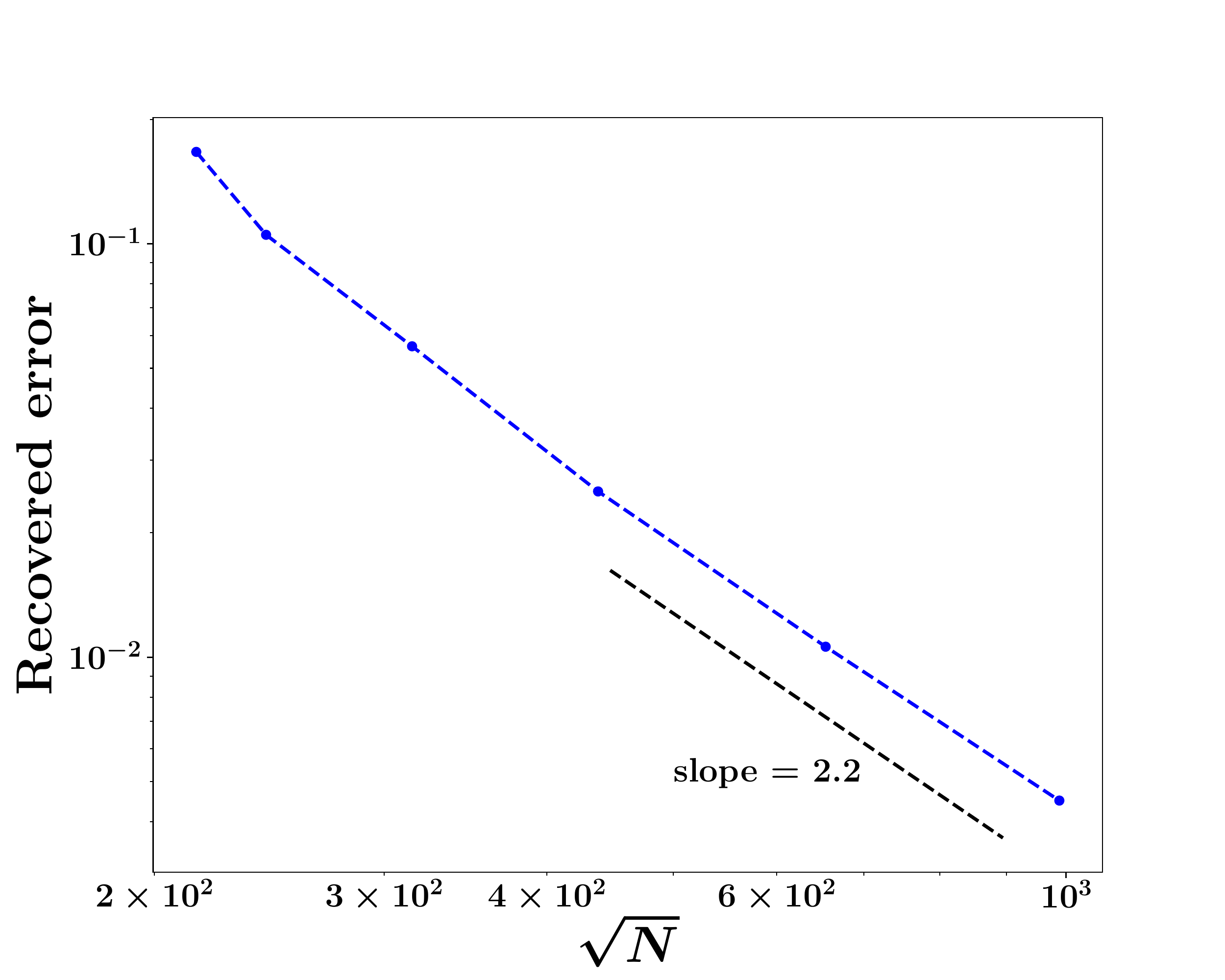}
        \caption{Convergence of the recovered error at the last time step. The slope measured by the last four data points is 2.2. $N$ denotes the total number of GMLS nodes.}
        \label{fig:100_circular_convergence}
    \end{subfigure}
    \quad
    \begin{subfigure}{.47\textwidth}
        \includegraphics[height=6cm]{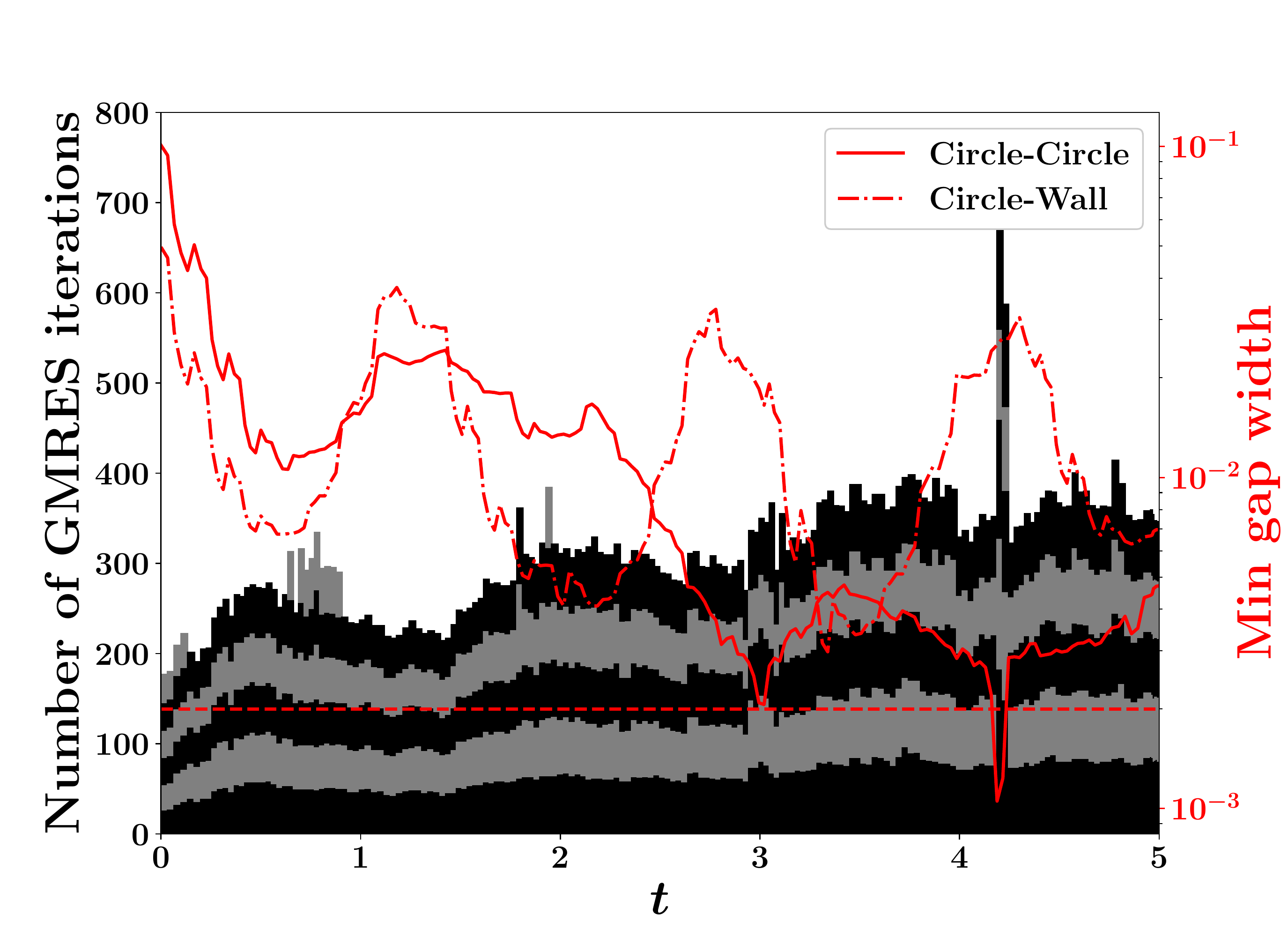}
        \caption{Number of GMRES iterations required in each adaptive $h$-refinement iteration and at each time step during the entire simulation.}
        \label{fig:100_circular_ite}
    \end{subfigure}
    \caption{Fluid-solid interactions$-$Particulate suspensions$-$100 similar particles: Convergence of the recovered error and the required number of GMRES iterations.}
\end{figure}

In the first simulation, the 100 particles are all circular with the radius $R = 0.04$. The snapshot of the particles' configuration at the terminal time $T = 5$ is shown in \Cref{fig:initial_configuration_100_circular}. The zoom-in pressure distributions at several locations are also shown in \Cref{fig:initial_configuration_100_circular}, for which we intentionally choose to show where the particles are in close contact either with each other or with the outer wall.
We can see that even though it is challenging in a dynamic simulation to resolve all point singularities governing lubrication effects, our numerical solver can stably predict the pressure fields without invoking any artificial subgrid-scale lubrication models. The convergence with respect to the recovered error, defined in Eq. \eqref{eq:global_posteriori_err}, for the last time step is shown in \Cref{fig:100_circular_convergence}. We find that the convergence rate reaches the theoretically expected 2nd order. 

To thoroughly examine the performance of our proposed preconditioner, we track the number of GMRES iterations required in each adaptive $h$-refinement iteration and in each time step, as depicted in \Cref{fig:100_circular_ite}. 
Several iterations of adaptive $h$-refinement are needed at each step to reach the preset error tolerance. Hence, the number of GMRES iterations required in each iteration of adaptive $h$-refinement at different time steps is rendered as the height of the bar with different colors and stacked together for various $h$-refinement iterations. For example, the lowest black bar represents the number of GMRES iterations required in the first iteration ($I=1$) of adaptive $h$-refinement at different time steps; the upper gray bar depicts the number of GMRES iterations needed for the second iteration ($I=2$) of adaptive $h$-refinement at different time steps; and so forth. By comparing the height of each bar at a fixed time step, we can see that the number of GMRES iterations generally stays constant across different iterations of adaptive $h$-refinement, indicating the scalability of the proposed monolithic GMG preconditioner in terms of increasing total DOFs. By comparing the heights of a bar at different time steps, we note that the number of GMRES iterations is highly correlated with the minimum gap width between solid boundaries (particle-particle or particle-wall). To elaborate on that, we add two lines (red) in \Cref{fig:100_circular_ite} showing the minimum gap widths between solid boundaries at different time instances. Note that the minimum gap width can be as small as $0.03R$ with $R = 0.04$ the particles' radius. 
Generally, when the minimum gap widths between solid boundaries are small, more GMRES iterations are required for the linear solver to converge. That is because finer GMLS nodes are needed to resolve narrower gaps between solid boundaries, resulting in more ill-conditioned linear systems to solve. The numbers of required GMRES iterations are comparable when the minimum gap widths are larger than 0.002 (indicated by the horizontal red dash line in \Cref{fig:100_circular_ite}), which is equivalent to $0.05R$ .

\subsubsection*{100 dissimilar particles}
\begin{figure}[htp]
    \centering
    \includegraphics[width=.98\textwidth]{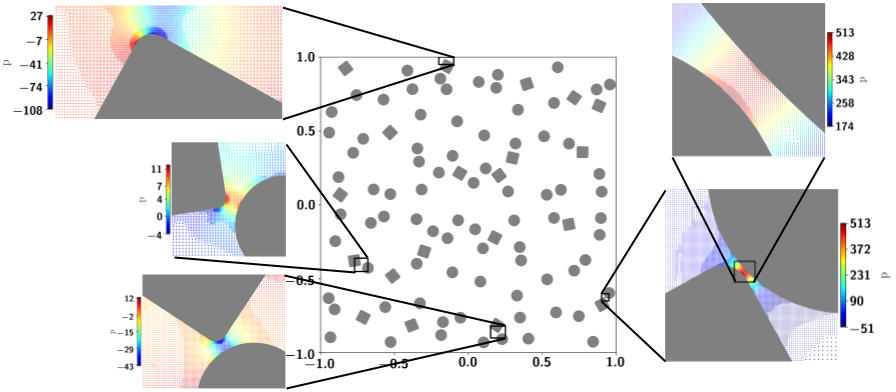}
    \caption{Fluid-solid interactions$-$Particulate suspensions: Configuration of 100 freely moving dissimilar particles in a Taylor–Green vortex flow at the terminal time. The zoom-in images are the computed pressure distributions at selected locations when the particles are either in close contact with each other or with the outer wall, where the color is correlated to the magnitude of pressure, and the point clouds are the GMLS nodes with adaptive refinement.}
    \label{fig:100_mix_configuration_contour}
\end{figure}

\begin{figure}[htp]
    \centering
    \begin{subfigure}{.47\textwidth}
    \centering
    \includegraphics[height=6cm]{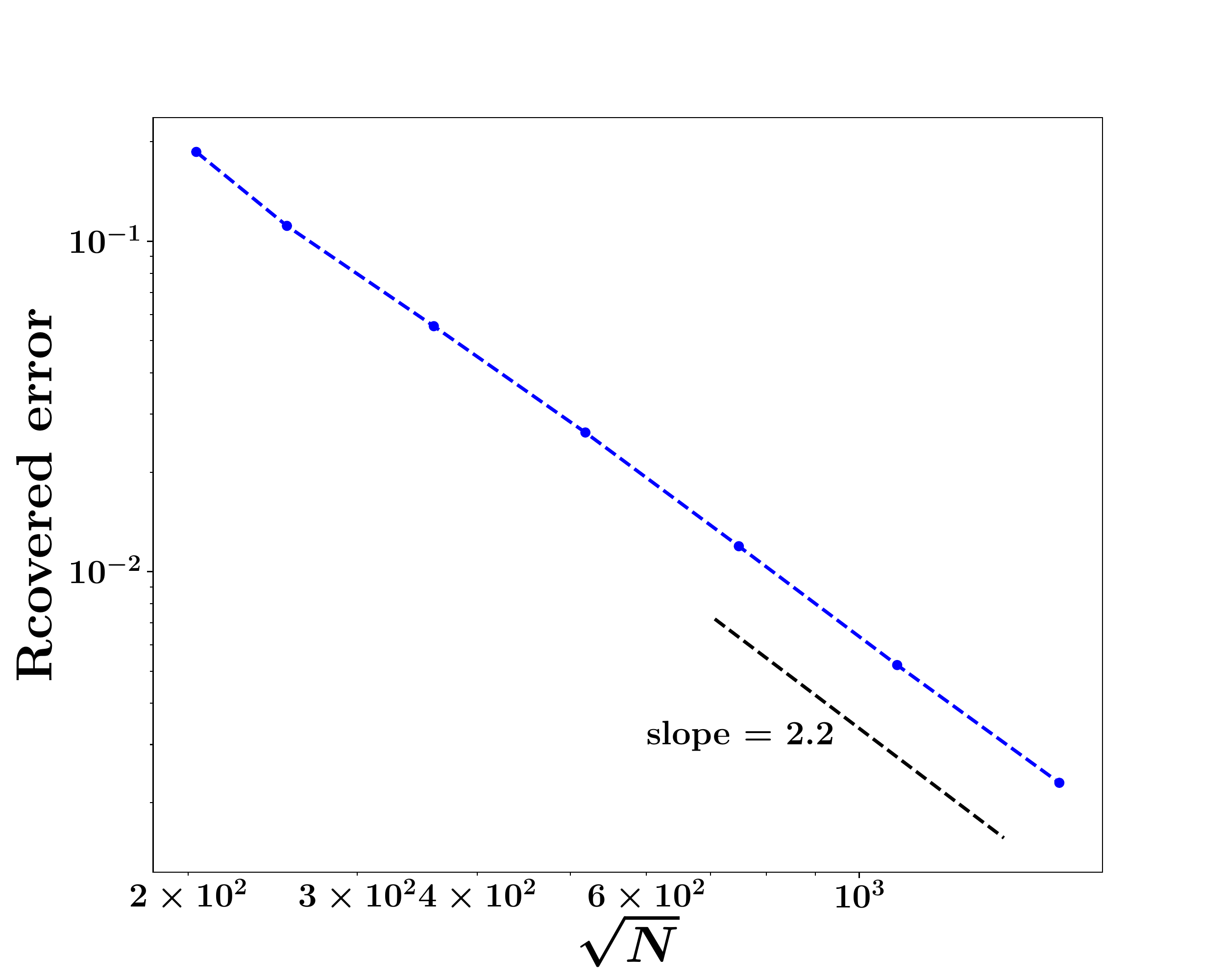}
    \caption{Convergence of the recovered error at the last time step. The slope measured by the last four data points is 2.2. $N$ denotes the total number of GMLS nodes.}
    \label{fig:100_mix_convergence}
    \end{subfigure}
    \quad \quad
    \begin{subfigure}{.47\textwidth}
    \centering
    \includegraphics[height=6cm]{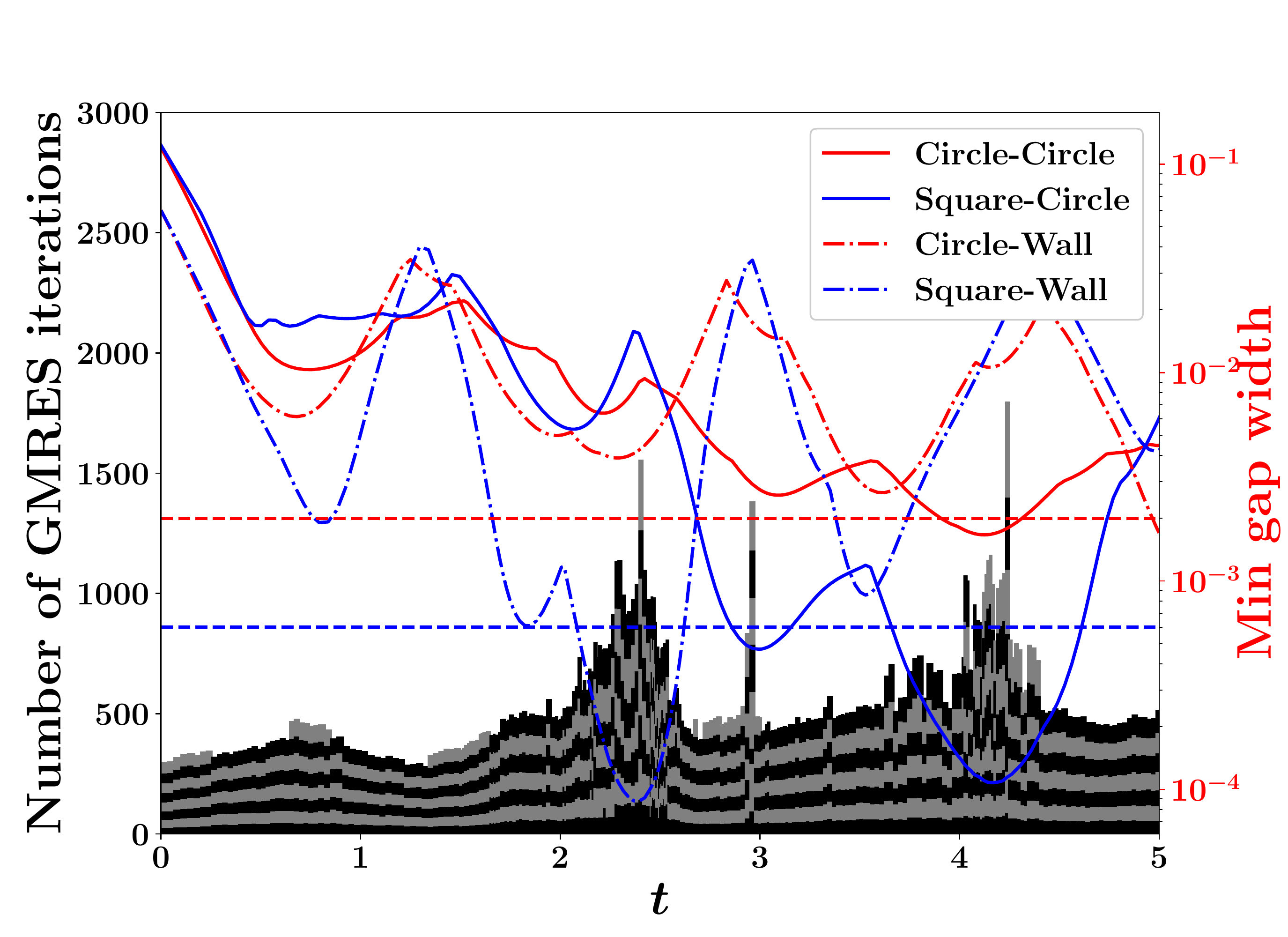}
    \caption{Number of GMRES iterations required in each adaptive $h$-refinement iteration and at each time step during the entire simulation.}
    \label{fig:100_mix_ite}
    \end{subfigure}
    \caption{Fluid-solid interactions$-$Particulate suspensions$-$100 dissimilar particles: Convergence of the recovered error and the required number of GMRES iterations.}
\end{figure}

In the second simulation, the suspended particles are a mixture of eighty circles and twenty squares. The radius of circular particles is still $R = 0.04$; the side length of square particles is $L=2R=0.08$. To mimic the particles in real applications of particulate suspensions, the squares are rounded at the corners with a rounding radius  $R^\prime = 0.1 L = 0.008$.
The snapshot of the particles' configuration at the terminal time $T = 5$ is shown in \Cref{fig:100_mix_configuration_contour}, where the zoom-in pressure distributions, particularly around square particles, are also shown at selected locations.
Regardless of particle shapes, our numerical solver can stably solve the problem even when the particles are in close contact with each other or with the outer wall boundary. No any artificial subgrid-scale lubrication model is employed during the entire simulation. The convergence with respect to the recovered error for the last time step is shown in \Cref{fig:100_mix_convergence}. With the inclusion of different shapes of solids, we still see the convergence rate reaching the theoretically expected 2nd order.  

By tracking the number of GMRES iterations required in each adaptive $h$-refinement iteration and in each time step, we assess the performance of the proposed preconditioning method. Similarly, we compare the number of GMRES iterations, as depicted in \Cref{fig:100_mix_ite}. For this case with mixed shapes of solids, more iterations of adaptive $h$-refinement are needed to reach the preset error tolerance at each time step. By comparing the height of each bar at a fixed time step, we also see that the number of GMRES iterations generally stays constant across different iterations of adaptive $h$-refinement, implying the scalability of the proposed monolithic GMG preconditioner in terms of increasing total DOFs. As in the case of 100 similar particles, the number of GMRES iterations is correlated with the minimum gap width between solid boundaries (particle-particle or particle-wall). In addition, the inclusion of square-shaped particles can worsen the problem's conditioning in the continuous limit. To reflect the combined effects of the minimum gap width between solid boundaries and the shape of particles, in \Cref{fig:100_mix_ite} we add lines to track the minimum gap widths associated with particles of different shapes (circular or square). We find that the square particles make the dominant contributions to the required GMRES iterations. When the minimum gap widths associated with square particles are larger than $0.0006$ (indicated by the horizontal blue dash line in \Cref{fig:100_mix_ite}), which is equivalent to $0.075 R^\prime$ with $R^\prime = 0.008$ the square particles' corner rounding radius, the number of GMRES iterations required at a given $h$-refinement iteration is generally comparable across different time instances; when they are smaller than $0.0006$, significantly more GMRES iterations are required for the linear solver to converge, because the corresponding linear system can become severely ill-conditioned.

\section{Conclusions} \label{sec:conclusion}
We have presented a monolithic GMG preconditioner for solving fluid-solid interaction problems in Stokes limit. The linear systems of equations are generated from the spatially adaptive GMLS discretization, which was developed in our previous work \cite{hu2019spatially}. The GMLS discretization is meshless and can handle large displacements and rotations of solid bodies without the expensive cost of generating and managing meshes. It guarantees the same order of accuracy for both the velocity and pressure fields, and high-order accuracy is achievable by its polynomial reconstruction property and choosing appropriate polynomial bases. A staggered discretization approximates the div-grad operator to ensure stable solutions of Stokes equations. With adaptive $h$-refinement directed by an \textit{a posteriori} recovery-based error estimator, it can resolve the singularities governing the lubrication effects between solid bodies without invoking any artificial subgrid-scale lubrication models. With the proposed monolithic GMG preconditioner in this work, we are able to scale up the spatially adaptive GMLS discretization for solving larger-scale fluid-solid interaction problems and achieve scalability with increasing numbers of solid bodies and total DOFs, while preserving accuracy and stability. 

The preconditioner is composed of two main ingredients: the interpolation/restriction operators and the smoothers. While the interpolation/restriction operators transfer the velocity and pressure values between different resolution levels, the smoothers damp the high-frequency error components on each resolution level. These two ingredients and their interplay determine the performance and scalability of the preconditioner and, thereby, the linear solver. The hierarchical structure of the adaptively refined GMLS nodes has provided us with an appropriate geometric setting for constructing the interpolation/restriction operators. In particular, to construct the interpolation operator, we employ local GMLS approximations from the coarse-level nodes to the fine-level nodes. To build the restriction operator, we let the value of a variable on a parent GMLS node equal the average value of its child nodes generated in a new iteration of adaptive refinement. During the interpolation and restriction processes, the divergence-free property of velocity is preserved, which plays an essential role in designing efficient MG preconditioning methods for solving Stokes equations. For the smoothers, we build a smoother for the fluid domain and a smoother for the solid bodies and then integrate them following a multiplicative overlapping Schwarz method. We have used a node-wise block Gauss-Seidel smoother for the fluid domain to address the coupling between the two field variables: velocity and pressure. The smoother for the solid bodies handles each solid body separately: for each solid body, a submatrix is first assembled corresponding to the solid body and its neighboring interior GMLS nodes; a Schur complement approach is then employed to approximately invert the submatrix based on the block splitting between the fluid and solid DOFs.

Such a developed preconditioner is integrated with the Krylov iterative solver, GMRES, for solving the linear systems of equations generated from the GMLS discretization. We have leveraged PETSc~\cite{petsc-user-ref} for the parallel implementation of the proposed monolithic GMG preconditioner. In addition to all associated linear algebra operations handled by PETSc, we have carefully taken care of domain decomposition, neighbor search, data storage, and imposing inhomogeneous Neumann BC for pressure in parallel implementation in order to warrant the parallel scalability of our numerical solver. Through a series of numerical tests, including simulating pure fluid flows and fluid-solid interactions with the inclusion of different numbers and shapes of solid bodies, we have demonstrated the performance and parallel scalability of our proposed preconditioner. Specifically, as the number of solid bodies and total DOFs increases, the convergence of the Krylov iterative solver can be ensured. For a fixed number of solid bodies, our preconditioner scales nearly linearly with respect to the total DOFs. When the number of solid bodies increases, our preconditioner exhibits sublinear optimality with respect to the number of solid bodies. In addition, our parallel implementation has achieved weak scalability for both pure fluid flows and fluid-solid interactions. For the flow in an artificial vascular network, the numerical result shows consistency with the prediction of Amdahl's law, indicating strong scalability and efficiency of our parallel implementation.

The present work has focused on 2D implementation and solving 2D problems. However, the proposed monolithic GMG preconditioning method is generally applicable to 2D or 3D fluid-solid interaction problems. Thus, extending the present work to 3D will be the next step in our future work. Although our current parallelization is sufficient for solving 2D problems, we expect that a more efficient parallel implementation is necessary when we move to 3D. That would include a better data storage strategy, more advanced domain decomposition and neighbor searching techniques, and parallelization of the monolithic GMG preconditioner accelerated by a hybrid GPU-CPU implementation. Furthermore, we have noted that the performance of the recovery-based \textit{a posteriori} error estimator in adaptive refinement can be affected by the PDE solutions' low regularity caused by the computational domain's nonconvexity, as discussed in the case of an artificial vascular network. As a result, it can lead to highly ill-conditioned linear systems. Therefore, developing a more theoretically guaranteed and robust \textit{a posteriori} error estimator would be another subject in our future work.

\section*{Declaration of Competing Interest}

The authors declare that they have no known competing financial interests or personal relationships that could have appeared to influence the work reported in this paper.

\section*{Acknowledgement}

Z. Y. was supported by the University of Wisconsin - Madison Office of the Vice Chancellor for Research and Graduate Education with funding from the Wisconsin Alumni Research Foundation; W. P. was supported by the National Science Foundation under Grant No. CMMI-1761068; X. H. was partially supported by the National Science Foundation under grant CCF-1934553.

\bibliographystyle{elsarticle-num}
\biboptions{numbers,sort&compress}
\bibliography{scalable_solver.bib}

\end{document}